\documentclass[10pt,twoside,leqno]{amsart}
\usepackage{amssymb,amsbsy,amsmath,amsfonts,amssymb,amscd,epsfig,times,graphics}
\newcommand{\C}			{\mathbb{C}}

\newcommand{\K}			{\mathbb{K}}
\newcommand{\R}			{\mathbb{R}}
\newcommand{\N}			{\mathbb{N}}

\newtheorem{theorem}{Theorem\it}
\newtheorem{lemma}{Lemma\it}
\newtheorem{proposition}{Proposition\it}

\newtheorem{definition}{Definition\rm}

\newtheorem{example}{\it Example\/\rm}

\begin{document}

\title[Symmetries of partial differential equations]{
Symmetries of partial differential equations
}

\author{Herv\'e Gaussier and Jo\"el Merker}

\address{(Gaussier and Merker) 
CNRS, Universit\'e de Provence, LATP, UMR 6632, CMI, 
39 rue Joliot-Curie, 13453 Marseille Cedex 13, France}

\email{[gaussier,merker]@cmi.univ-mrs.fr} 


\date{\number\year-\number\month-\number\day}

\begin{abstract}
We establish a link between the study of completely integrable systems
of partial differential equations and the study of generic
submanifolds in $\C^n$. Using the recent developments of
Cauchy-Riemann geometry we provide the set of symmetries of such a
system with a Lie group structure. Finally we determine the precise
upper bound of the dimension of this Lie group for some specific
systems of partial differential equations.
\end{abstract}

\maketitle

\begin{center}
\begin{minipage}[t]{10cm}
\baselineskip =0.35cm
{\scriptsize

\centerline{\bf Table of contents}

\smallskip

{\bf 1.~Introduction \dotfill 1.}

{\bf 2.~Submanifold of solutions \dotfill 3.}

{\bf 3.~ Lie theory for partial differential equations\dotfill11.}

{\bf 4.~Optimal upper bound on $\dim_\K\mathfrak{Sym}(\mathcal{E})$
when $n=m=1$ \dotfill 17.}

{\bf 5.~Optimal upper bound on $\dim_\K\mathfrak{Sym}(\mathcal{E})$ in
the general dimensional case\dotfill 20.}

}\end{minipage}
\end{center}

\bigskip

\section{Introduction}
To study the geometry of a real analytic Levi nondegenerate
hypersurface $M$ in $\C^{2}$, one of the principal ideas of
H.~Poincar\'e, of B.~Segre and of \'E.~Cartan in the fundamental
memoirs \cite{po1907}, \cite{se1931}, \cite{se1932}, \cite{ca1932} was to
associate to $M$ a system $(\mathcal E_M)$ of (partial) differential
equations, in order to solve the so-called equivalence
problem. Establishing a natural correspondence between the local
holomorphic automorphisms of $M$ and the Lie symmetries of $(\mathcal
E_M)$ they could use the classification results on differential
equations achieved by S.~Lie in \cite{enlie1890} and pursued by
A.~Tresse in \cite{tr1896}.

Starting with such a correspondence, we shall establish a general link
between the study of a real analytic generic submanifold of
codimension $m$ in $\C^{n+m}$ and the study of completely integrable
systems of analytic partial differential equations. We shall observe
that the recent theories in Cauchy-Riemann (CR) geometry may be
transposed to the setting of partial differential equations, providing
some new information on their Lie symmetries. 

Indeed consider for $\K=\R$ or $\C$ a $\K$-analytic system
$(\mathcal E)$ of the following general form:

\vskip 0,1cm
\noindent ($\mathcal E$) \centerline{
$u_{x^\alpha}^j(x) =
F_\alpha ^j\left(x,u(x),(u_{x^{\beta(q)}}^{j(q)}(x))_{1 \leq q \leq p}\right).$
}
\vskip 0,1cm
\noindent Here $x = (x_1,\dots,x_n)\in \K^n$,
$u=(u^1,\dots,u^m)\in\K^m$, the integers $j(1),\dots, j(p)$ satisfy $1
\leq j(q) \leq m$ for $q=1,\dots,p$, and $\alpha$ and the multiindices
$\beta(1),\dots, \beta(p) \in \N^n$ satisfy $|\alpha|,\ |\beta(q)|
\geq 1$. We also require $(j,\alpha) \neq (j(1),\beta(1)),\dots,
(j(p),\beta(p))$. For $j= 1, \dots, m$ and $\alpha \in \N^n$, we
denote by $u_{ x^\alpha}^j$ the partial derivative $\partial^{ \vert
\alpha \vert} u^j/ \partial x^\alpha$. We assume that the system
$(\mathcal E)$ is {\sl completely integrable}, namely that the
Pfaffian system naturally associated in the jet space is involutive in
the sense of Frobenius. We note that in that case ($\mathcal E$) is
{\sl locally solvable}, meaning that through every point
$(x^*,u^*,u_\beta^*, u_\alpha^*)$ in the jet space, satisfying
$u_\alpha^*=F_\alpha(x^*,u^*,u_\beta^*)$ (written in a condensed
form), there exists a local $\K$-analytic solution $u=u(x)$ of
($\mathcal E$) satisfying $u(x^*) = u^*$ and
$u_{x^\beta}(x^*)=u_\beta^*$. Consequently the Lie theory
(\cite{ol1986}) may be applied to such systems.  We shall associate
with $(\mathcal E)$ the {\sl submanifold of solutions} $\mathcal M$ in
$\K^{n+2m+p}$ given by $\K$-analytic equations of the form
\begin{equation}\label{e12}
u^j=\Omega_j(x,\nu,\chi),\ \ j=1,\dots,m,
\end{equation}
where $\nu \in \K^m$ and where $\chi \in \K^p$. Moreover the integer
$m+p$ is the number of initial conditions for the general solution
$u(x):=\Omega(x,\nu,\chi)$ of $(\mathcal E)$, whose existence and
uniqueness follow from complete integrability. Precisely, the
parameters $\nu,\chi$ correspond to the data $u(0),
(u_{x^{\beta(q)}}^{j(q)}(0))_{1 \leq q \leq p}$. In the special case
where the system $(\mathcal E)$ is constructed from a generic
submanifold $M$ as in \cite{se1931}, \cite{su2001a} (see also
Subsection~2.2 below), the corresponding submanifold of solutions is
exactly the extrinsic complexification of $M$.

A pointwise $\K$-analytic transformation $(x',u')=\Phi(x,u)$ defined
in a neighbourhood of the origin and sufficiently close to the
iedntity mapping is called a {\sl Lie symmetry of} $(\mathcal E)$ if
it transforms the graph of every solution to the graph of an other
local solution. A vector field $X=\sum_{l=1}^n \, Q^l(x,u)\, \partial
/\partial x_l+ \sum_{j=1}^m \, R^j(x,u)\, \partial /\partial u^j$ is
called an {\sl infinitesimal symmetry of} $(\mathcal{E})$ if for every
$s$ close to zero in $\K$ the local diffeomorphism $(x,u)\mapsto
\exp(sX)(x,u)$ associated to the flow of $X$ is a Lie symmetry of
$\mathcal{E}$.  According to~\cite{ol1986} (Chapter~2) the
infinitesimal symmetries of $(\mathcal E)$ form a Lie algebra of
vector fields defined in a neighbourhood of the origin in $\K^n\times
\K^m$, denoted by $\mathfrak{Sym}(\mathcal{E})$. Inspired by recent
developments in CR geometry we shall provide in Section~2
nondegeneracy conditions on $\mathcal M$ insuring firstly that
$\mathfrak{Sym}(\mathcal{E})$ may be identified with the Lie algebra
$\mathfrak{Sym}(\mathcal M)$ of vector fields of the form
\begin{equation}\label{eee20}
\sum_{l=1}^n\, 
Q^l(x,u)\, \frac{\partial}{\partial x_l}+\sum_{j=1}^m\,
R^j(x,u)\, \frac{\partial}{\partial u^j}+
\sum_{j=1}^m\, 
\Pi^j(\nu,\chi)\,
\frac{\partial }{\partial \nu^j}+
\sum_{q=1}^p\, 
\Lambda^q(\nu,\chi)\, 
\frac{\partial }{\partial \chi_q},
\end{equation}
which are tangent to $\mathcal{M}$, and secondly that
$\mathfrak{Sym}(\mathcal M) \cong \mathfrak{Sym}(\mathcal E)$ is
finite dimensional. The strength of this identification is to provide
some (non optimal) bound on the dimension of $\mathfrak{Sym}(\mathcal
E)$ for arbitrary systems of partial differential equations with an
arbitrary number of variables, see~Theorem~\ref{tthm23}.

In the second part of the paper (Sections~3, 4 and 5), using the
classical Lie theory ({\it cf}. \cite{enlie1890}, \cite{ol1986},
\cite{ol1995} and \cite{bk1989}), we provide an optimal upper bound on
the dimension of $\mathfrak{Sym}(\mathcal{E})$ for a completely
integrable $\K$-analytic system $(\mathcal E)$ of the following form:
\vskip 0,1cm

\noindent ($\mathcal{E}$) \centerline{
$u_{x^\alpha}^j = F_\alpha^j(x,u(x),(u_{x^\beta}(x))_{1 \leq |\beta|
\leq \kappa-1}), \ \ \ \ \ \alpha \in \N^n,\ |\alpha| = \kappa,\
j=1,\dots,m$.
}
\vskip 0,1cm
\noindent This system is a special case of the system studied in
Section~2.  For instance the {\sl homogeneous} system $(\mathcal E_0)
\ :\ \ u_{x_{k_1}\cdots x_{k_\kappa}}^j(x)=0$ is completely
integrable. The solutions of ($\mathcal E_0$) are the polynomials of
the form $u^j(x)=\sum_{ \beta \in \N^n,\, \vert \beta \vert \leq
\kappa-1}\, \lambda_\beta^j\, x^\beta$, $j=1, \dots, m$, where
$\lambda_\beta^j \in\K$ and a Lie symmetry of $(\mathcal E_0)$ is a
transformation stabilizing the graphs of polynomials of degree $\leq
\kappa-1$.  We prove the following Theorem:

\begin{theorem}\label{thm1}
Let $(\mathcal E)$ be the $\K$-analytic system of partial differential
equations of order $\kappa \geq 2$, with $n$ independent variables and
$m$ dependent variables, defined just above. Assume that $(\mathcal
E)$ is completely integrable. Then the Lie algebra
$\mathfrak{Sym}(\mathcal{E})$ of its infinitesimal symmetries
satisfies the following estimates:
\begin{equation}\label{e13}
\left\{\begin{array}{llll}
\dim_\K(\mathfrak{Sym}(\mathcal E)) &\leq & (n+m+2)(n+m), 
& \ \ {\rm if}\ \ \kappa = 2,
\\
\dim_\K(\mathfrak{Sym}(\mathcal E)) &\leq & 
n^2+2n+m^2+
m \, C_{n+\kappa -1}^{\kappa -1}, & \ \ {\rm if}\ \ \kappa \geq 3,
\end{array}
\right.
\end{equation}
where we denote $C_{n+\kappa-1}^{\kappa-1}:= \frac{(n+\kappa-1)!}{n! \
(\kappa-1)!}$. Moreover the inequalities~(\ref{e13}) become equalities
for the homogeneous system $(\mathcal{E}_0)$.
\end{theorem}
\noindent We remark that there is no combinatorial formula
interpolating these two estimates. Theorem~\ref{thm1} is a
generalization of the following results. For $n=m=1$, S.~Lie proved
that the dimension of the {\sl Lie} algebra $\mathfrak{Sym}
\,(\mathcal{ E})$ is less than or equal to 8 if $\kappa =2$ and is
less than or equal to $\kappa+4$ if $\kappa \geq 3$, these bounds
being reached for the homogeneous system~({\it cf.}~\cite{enlie1890}).
For $n=1$, $m\geq 1$ and $\kappa=2$, F.~Gonz\'alez-Gasc\'on and
A.~Gonz\'alez-L\'opez proved in~\cite{gg1983} that the dimension of
$\mathfrak{Sym}\, (\mathcal{ E})$ is less than or equal to $(m+3)
(m+1)$.  For $n=1$, $m\geq 1$ and $\kappa=2$, using the equivalence
method due to \'E. Cartan, M.~Fels~\cite{fe1995} proved that the
dimension of $\mathfrak{Sym} \, (\mathcal{E})$ is less than or equal
to $m^2+4m+3$, with equality if and only if the system ($\mathcal E$)
is equivalent to the system $u_{x^2}^j=0$, $j=1,\dots,m$.  He also
generalized this result to the case $n=1$, $m\geq 1$, $\kappa=3$.  For
$n\geq 1$, $m\geq 1$ and $\kappa =2$, A.~Sukhov proved
in~\cite{su2001a} that the dimension of $\mathfrak{Sym} \,
(\mathcal{E})$ is less than or equal to $(n+m+2) (n+m)$ (the first
inequality in Theorem~\ref{thm1}), with equality for the homogeneous
system $u_{x_{k_1}x_{k_2}}^j=0$. 

Consequently, for the case $\kappa =2$, we will only give the general
form of the Lie symmetries of the homogeneous system ($\mathcal E_0$)
(see Subsection~5.2). We will prove Theorem~\ref{thm1} for the case
$\kappa \geq 3$. The formulas obtained in Sections~3, 4 and 5 were
checked with the help of MAPLE release 6.

\vskip 0,1cm
\noindent{\sl Acknowledgment.}  This article was written while the
first author had a six months delegation position at the CNRS. He
thanks this institution for providing him this research
opportunity. The authors are indebted to G\'erard Henry, the computer
ing\'enieur (LATP, UMR 6632 CNRS), for his technical support.

\section{Submanifold of solutions}

\subsection{Preliminary} 
Let $\K=\R$ or $\C$.  Let $n\geq 1$ and let $x= (x_1, \dots, x_n) \in
\N$. We denote by $\K\{x\}$ the local ring of $\K$-analytic functions
$\varphi=\varphi(x)$ defined in some neighbourhood of the origin in
$\K^n$.  If $\varphi \in \K\{x\}$ we denote by $\bar \varphi$ the
function in $\K\{x\}$ satisfying $\overline{\varphi(x)}\equiv \bar
\varphi(\bar x)$. Recall that a $\K$-analytic function $\varphi$
defined in a domain $U\subset \K^n$ is called {\sl $\K$-algebraic} (in
the sense of Nash) if there exists a nonzero polynomial $P= P( X_1,
\dots, X_n, \Phi) \in \K [X_1, \dots, X_n, \Phi]$ such that $P (x,
\varphi (x)) \equiv 0$ on $U$. All the considerations in this paper
will be local: functions, submanifolds and mappings will always be
defined in a small connected neighbourhood of some point (most often
the origin) in $\K^n$.

\subsection{System of partial differential equations
associated to a generic submanifold of $\C^{n+m}$} Let $M$ be a real
algebraic or analytic local submanifold of codimension $m$ in
$\C^{n+m}$, passing through the origin. We assume that $M$ is {\sl
generic}, namely $T_0M+iT_0M= T_0\C^{n+m}$. Classically ({\it
cf.}~\cite{ber1999}) there exists a choice of complex linear
coordinates $t=(z,w)\in\C^n\times\C^m$ centered at the origin such
that $T_0M=\{{\rm Im}\, w=0\}$ and such that there exist $m$ complex
algebraic or analytic defining equations representing $M$ as the set
of $(z,w)$ in a neighbourhood of the origin in $\C^{n+m}$ which
satisfy
\begin{equation}\label{eee21}
w_1=\Theta_1(z,\bar z,\bar w), \dots\dots,
w_m=\Theta_m(z,\bar z,\bar w).
\end{equation}
Furthermore, the mapping $\Theta=(\Theta_1,\dots,\Theta_m)$ satisfies
the functional equation
\begin{equation}\label{eee22}
w\equiv \Theta(z, \bar z, \overline{\Theta}(\bar z, z, w)),
\end{equation}
which reflects the reality of the generic submanifold $M$. It follows
in particular from~(\ref{eee22}) that the local holomorphic mapping
$\C^m\ni \bar w \mapsto (\Theta_j (0, 0, \bar w))_{1 \leq j \leq m}
\in \C^m$ is of rank $m$ at $\bar w=0$.

Generalizing an idea due to B.~Segre in~\cite{se1931} and
\cite{se1932}, exploited by \'E.~Cartan in~\cite{ca1932} and
more recently by A.~Sukhov in~\cite{su2001a},
\cite{su2001b}, \cite{su2002}, we shall associate to $M$ a system of
partial differential equations. For this, we need some
general nondegeneracy condition, which generalizes Levi
nondegeneracy. Let $\ell_0\in\N$ with $\ell_0\geq 1$. We shall assume
that $M$ is $\ell_0$-finitely nondegenerate at the origin, {\it
cf.}~\cite{ber1999}, \cite{me2003}, \cite{gm2002a}. This means that
there exist multiindices $\beta(1), \dots, \beta(n) \in \N^n$ with
$\vert \beta(k) \vert \geq 1$ for $k=1, \dots, n$ and $\max_{1\leq k
\leq n}\, \vert \beta (k) \vert =\ell_0$, and integers
$j(1),\dots,j(n)$ with $1 \leq j(k) \leq m$ for $k=1, \dots, n$ such
that the local holomorphic mapping
\begin{equation}\label{eee23}
\C^{n+m}\ni(\bar z, \bar w) \longmapsto \left( 
(\Theta_j(0,\bar z, \bar w))_{
1\leq j\leq m}, 
\left(
\Theta_{j(k),z^{\beta(k)}}(0,\bar z,\bar w)\right)_{1\leq 
k\leq n}
\right)\in\C^{m+n}
\end{equation}
is of rank equal to $n+m$ at $(\bar z,\bar w)=(0,0)$.  Here, we denote
the partial derivative $\partial^{ \vert \beta \vert} \Theta_j (0,
\bar z, \bar w) / \partial z^\beta$ simply by $\Theta_{ j, z^\beta }(
0, \bar z, \bar w )$. Then $M$ is Levi nondegenerate at the origin if
and only if $\ell_0=1$. By complexifying the variables $\bar
z$ and $\bar w$, we get new independent variables $\zeta\in\C^n$ and
$\xi\in\C^m$ together with a complex algebraic or analytic
$m$-codimensional submanifold $\mathcal{M}$ in $\C^{2(n+m)}$ of
equations
\begin{equation}\label{eee24}
w_j=\Theta_j(z,\zeta,\xi), \ \ \ \ \ \ \ \ \
j=1,\dots,m,
\end{equation}
called the {\sl extrinsic complexification of $M$}.  In
the defining equations~(\ref{eee24}) of $\mathcal{ M}$,
following~\cite{se1931} and~\cite{su2001a}, we may consider the
``dependent variables'' $w_1, \dots, w_m$ as algebraic or analytic
functions of the ``independent variables'' $z=(z_1, \dots, z_n)$, with
additional dependence on the extra ``parameters''
$(\zeta,\xi)\in\C^{n+m}$. Then by applying the differential operator
$\partial^{\vert \alpha \vert} /\partial z^\alpha$ to~(\ref{eee24}), we
obtain $w_{j,z^\alpha}(z)= \Theta_{j,z^\alpha}(z,\zeta,\xi)$. Writing
these equations for $(j,\alpha)=(j(k),\beta(k))$ with $k=1,\dots,n$, we
obtain a system of $m+n$ equations
\begin{equation}\label{eee25}
\left\{
\aligned
w_j(z)=
& \
\Theta_j(z,\zeta,\xi),
\ \ \ \ \ j=1,\dots,m,\\
w_{j(k),z^{\beta(k)}}(z)=
& \
\Theta_{j(k),z^{\beta(k)}}
(z,\zeta,\xi), \ \ \ \ \ k=1,\dots,n.
\endaligned\right.
\end{equation}
In this system~(\ref{eee25}), by the assumption of $\ell_0$-finite
nondegeneracy~(\ref{eee23}), the algebraic or analytic implicit
function theorem allows to solve the parameters $(\zeta, \xi)$ in
terms of the variables $(z_k, w_j(z), w_{j(k), z^{\beta (k)}} (z))$,
providing a local algebraic or analytic $\C^{ n+m}$-valued mapping $R$
such that $(\zeta, \xi)= R\left(z_k, w_j(z), w_{j(k), z^{\beta (k)}}
(z)\right)$.  Finally, for every pair $(j,\alpha)$ different from
$(1,0), \dots, (m,0), (j(1), \beta(1)), \dots, (j(n), \beta(n))$, we
may replace $( \zeta, \xi)$ by $R$ in the differentiated expression
$w_{ j, z^\alpha} (z) = \Theta_{ j, z^\alpha} (z, \zeta, \xi)$. This
yields
\begin{equation}\label{eee26}
\aligned
w_{j,z^\alpha}(z)
& \
=
\Theta_{j,z^\alpha}\left(z, R(z_k,w_j(z),w_{j(k),z^{\beta(k)}}(z))\right)\\
& \
=: 
F_{j,\alpha}\left(z_k,w_j(z),w_{j(k),z^{\beta(k)}}(z)\right).
\endaligned
\end{equation}
This is the {\sl system of partial differential equations associated
with $\mathcal{M}$}. As argued by B.~Segre in~\cite{se1931},
the geometric study of generic submanifolds of
$\C^n$ may gain much information from the study of their associated
systems of partial differential equations ({\it cf.}~\cite{su2001a},
\cite{su2001b}).  The next paragraphs are devoted to provide a {\it
general one-to-one correspondence}\, between completely integrable
systems of analytic partial differential equations and their
associated ``submanifolds of solutions'' (to be defined precisely
below) like $\mathcal{M}$ above. Afterwards, we shall observe that
conversely, the study of systems of analytic partial differential
equations also gains much information from the direct study of their
associated submanifolds of solutions.

\subsection{Completely integrable systems of partial differential 
equations} Let now $n, \ m,\ p \in \N$ with $n,\ m,\ p \geq 1$, let
$\kappa \in \N$ with $\kappa \geq 2$ and let $u= (u^1, \dots, u^m) \in
\K^m$. Consider a collection of $p$ multiindices $\beta (1), \dots,
\beta(p) \in \N^n$ with $\vert \beta(q) \vert \geq 1$ for
$q=1,\dots,p$ and $\max_{1\leq q \leq p}\, \vert \beta(q) \vert
=\kappa-1$. Consider also $p$ integers $j(1), \dots, j(p)$ with $1
\leq j(q) \leq m$ for $q= 1, \dots, p$. Inspired by~(\ref{eee26}), we
consider a general system of partial differential equations of $n$
independent variables $(x_1,\dots,x_n)$ and $m$ dependent variables
$(u^1,\dots,u^m)$ which is of the following form:
\vskip 0,1cm
\noindent ($\mathcal{E}$) \centerline{
$u_{x^\alpha}^j(x)=F_\alpha^j\left(x,
u(x),(u_{x^{\beta(q)}}^{j(q)}(x))_{1\leq 
q\leq p}\right),$
}

\vskip 0,1cm
\noindent where $(j,\alpha)\neq (j(1),\beta(1)),\dots,(j(p),\beta(p))$
and $j=1,\dots,m$, $\vert \alpha \vert \leq \kappa$. Here, we assume
that $u=0$ is a local solution of the system $(\mathcal{E})$ and that
the functions $F_\alpha^j$ are $\K$-algebraic or $\K$-analytic in a
neighbourhood of the origin in $\K^{n+m+p}$.  Among such systems are
included ordinary differential equations of any order $\kappa\geq 2$,
systems of second order partial differential equation as studied
in~\cite{su2001a}, {\it etc.}

Throughout this article, we shall assume the system $(\mathcal{E})$
completely integrable. By analyzing the application of the Frobenius
theorem in jet spaces, one can show (we will not develop this) that
the general solution of the system ($\mathcal{E}$) is given by
$u(x):=\Omega(x,\nu,\chi)$, where the parameters $\nu \in \K^n$ and
$\chi \in \K^n$ essentially correspond to the ``initial conditions''
$u(0)$ and $(u_{x^{\beta(q)}}^{j(q)}(0))_{1\leq q\leq p}$, and
$\Omega$ is a $\K$-analytic $\K^n$-valued mapping. In the case of a
generic submanifold as in Subsection~2.2 above, we recover the mapping
$\Theta$. In the sequel, we shall use the following terminology: the
coordinates $(x,u)$ will be called the {\sl variables} and the
coordinates $(\nu,\chi)$ will be called the {\sl parameters} or the
{\sl initial conditions}. In Subsection~2.5 below, we shall introduce
a certain duality where the r\^oles between variables and parameters
are exchanged.

\subsection{Associated submanifold of solutions}
The existence of the function $\Omega$ and the analogy with Subsection~2.2
leads us to introduce the {\sl submanifold of solutions associated to
the completely integrable system $(\mathcal{E})$}, which by definition
is the $m$-codimensional $\K$-analytic submanifold of $\K^{n+2m+p}$,
equipped with the coordinates $(x,u,\nu,\chi)$, defined by the
Cartesian equations
\begin{equation}\label{eee27}
u_j=\Omega_j(x,\nu,\chi), \ \ \ \ \ \ \ \ 
j=1,\dots,m.
\end{equation}
Let us denote this submanifold by $\mathcal{M}$.  We stress that
in general such a submanifold cannot coincide with the
complexification of a generic submanifold of $\C^{m+n}$, for instance
because $\K$ may be equal to $\R$ or, if $\K=\C$, because the integer
$p$ is not necessarily equal to $n$. Also, even if $\K=\C$ and $n=p$,
the mapping $\Omega$ does not satisfy a functional equation
like~(\ref{eee22}). In fact, it may be easily established that the
submanifold of solutions of a completely integrable system
of partial differential equations like $(\mathcal{E})$ coincides with
the complexification of a generic submanifold {\it if and only 
if}\, $\K=\C$, $p=n$ and the mapping $\Omega$ satisfies
a functional equation like~(\ref{eee22}).

Let now $\mathcal{M}$ be a
submanifold of $\K^{n+2n+p}$ of the form~(\ref{eee27}), but not
necessarily constructed as the submanifold of solutions of a system
$(\mathcal{E})$. We shall always assume that
$\Omega_j(0,\nu,\chi)\equiv \nu^j$.  We say that $\mathcal{M}$ is {\sl
solvable with respect to the parameters} if there exist
multiindices $\beta(1),\dots,\beta(p)\in\N^n$ with $\vert \beta(q)
\vert \geq 1$ for $q=1,\dots,p$ and integers $j(1),\dots,j(p)$ with
$1\leq j(q) \leq m$ for $q=1,\dots,p$ such that the local
$\K$-analytic mapping
\begin{equation}\label{eee28}
\K^{m+p}\ni (\nu,\chi)\longmapsto 
\left(
(\Omega_j(0,\nu,\chi)_{1\leq j\leq m}, \ 
\left(
\Omega_{j(q), x^{\beta(q)}}(0,\nu,\chi)
\right)_{1\leq q \leq p}
\right)\in\K^{m+p}
\end{equation}
is of rank equal to $m+p$ at $(\zeta,\chi)=(0,0)$ (notice that since
$\Omega_j(0,\nu,\chi)\equiv \nu^j$, then the $m$ first components of
the mapping~(\ref{eee28}) are already of rank $m$). We remark that the
submanifold of solutions of a system $(\mathcal{E})$ is automatically
solvable with respect to the variables, the multiindices $\beta(q)$
and the integers $j(q)$ being the same as in the arguments of the
right hand side terms $F_\alpha^j$ in $(\mathcal{E})$.

\subsection{Dual system of defining equations}
Since $\Omega_j(0,\nu,\chi)\equiv \nu^j$, we may solve the
equations~(\ref{eee27}) with respect to $\nu$ by means of the analytic
implicit function theorem, getting an equivalent system of equations
for $\mathcal{M}$:
\begin{equation}\label{eee29}
\nu^j=\Omega_j^*(\chi,x,u), \ \ \ \ \ \ \ \ \  \ j=1,\dots,m.
\end{equation}
We call this the {\sl dual system of defining equations for $\mathcal{M}$}.
By construction, we have the functional equation
\begin{equation}\label{eee210}
u\equiv 
\Omega(x,\Omega^*(\chi,x,u),\chi),
\end{equation}
implying the identity $\Omega_j^*(0,x,u)\equiv u^j$.
We say that $\mathcal{M}$ is {\sl solvable with respect to the variables}
if there exist multiindices $\delta(1),\dots,\delta(n)\in\N^p$
with $\vert \delta(l)\vert \geq 1$ for $l=1,\dots,n$ and integers
$j(1),\dots,j(n)$ with $1\leq j(l) \leq m$ for $l=1,\dots,m$ such that
the local $\K$-analytic mapping
\begin{equation}\label{eee211}
\K^{n+m} \ni (x,u) \longmapsto 
\left(
(\Omega_j^*(0,x,u))_{1\leq j\leq m}, \
\left(
\Omega_{j(l), \, \chi^{\delta(l)}}^*(0,x,u)
\right)_{1\leq l\leq n}
\right)\in\K^{m+n}
\end{equation}
is of rank equal to $n+m$ at $(x,u)=(0,0)$ (notice that since
$\Omega_j^*(0,x,u)\equiv u^j$, the $m$ fisrt components of the
mapping~(\ref{eee211}) are already of rank $m$).

In the case where $\mathcal M$ is the complexification of a generic
submanifold then the solvability with respect to the parameters is
equivalent to the solvability with respect to the variables since
$\Omega^* \equiv \overline{\Omega}$. However we notice that a
submanifold $\mathcal{M}$ of solutions of a system $(\mathcal{E})$ is
not automatically solvable with respect to the variables, as shows the
following trivial example.

\begin{example}\label{ex21}
{\rm 
Let $n=2$, $m=1$ and let $(\mathcal{E})$ denote the system $u_{x_2}=0$,
$u_{x_1x_1}=0$, whose general solutions are $u(x)=\nu+x_1\chi=:
\Omega(x_1,x_2,\nu,\chi)$. Notice that the variable $x_2$ is absent
from the dual equation $\nu=u-x_1\chi_1=:\Omega^*(\chi,x_1,x_2,u)$.
It follows that $\mathcal{M}$ is not solvable with respect to the
variables.  }
\end{example}
 
\subsection{Symmetries of $(\mathcal{E})$, their
lift to the jet space and their lift to the parameter space} We denote
by $\mathcal{J}_{n,m}^\kappa$ the space of jets of order $\kappa$ of
$\K$-analytic mappings $u=u(x)$ from $\K^n$ to $\K^m$.  Let
\begin{equation}\label{eee212}
(x_l,u^j,U_{l_1}^{i_1}, 
U_{l_1,l_2}^{i_1},\dots,
U_{l_1,\dots,l_\kappa}^{i_1})\in \K^{n+m\, C_{\kappa+n}^{\kappa}}
\end{equation}
denote the natural coordinates on $\mathcal{J}_{ n,m}^\kappa$. Here,
the superscripts $j,i_1$ and the subscripts $l, l_1, l_2, \dots,
l_\kappa$ satisfy $j, i_1=1,\dots,m$ and $l, l_1,
l_2, \dots, l_\kappa=1, \dots,n$. The independent coordinate $U_{ l_1,
\dots, l_\lambda}^{i_1}$ corresponds to the partial derivative $u_{
x_{l_1} \dots x_{ l_\lambda }}^{ i_1}$. Finally, by symmetry of
partial differentiation, we identity every coordinate $U_{l_1, \dots,
l_\lambda}^{ i_1}$ with the coordinates $U_{\sigma( l_1), \dots,
\sigma( l_\lambda)}^{ i_1}$, where $\sigma$ is an arbitrary
permutation of the set $\{1, \dots, \lambda\}$.  With these
identifications, the $\kappa$-th order jet space $\mathcal{ J}_{
n,m}^\kappa$ is of dimension $n+m \, C_{ \kappa+n }^\kappa$, where
$C_p^q :=\frac{ p!}{ q! \ (p-q)!}$ denotes the binomial
coefficient. Also, we shall sometimes use an equivalent notation for
coordinates on $\mathcal{ J}_{ n,m}^\kappa$:
\begin{equation}\label{eee213}
(x_l,u^j, U_\beta^i)\in\K^{n+m\, C_{\kappa+n}^n},
\end{equation} 
where $\beta\in\N^n$ satisfies $\vert \beta \vert \leq \kappa$ and
where the independent coordinate $U_\beta^i$ corresponds to the
partial derivative $u_{x^\beta}^i$.   

associated to the system $(\mathcal{ E})$ is the so-called {\sl skeleton}
$\Delta_\mathcal{ E}$, which is the $\K$-analytic submanifold
of dimension $n+m+p$ in
$\mathcal{J}_{n,m}^\kappa$ simply defined by replacing the partial
derivatives of the dependent
variables $u^j$ by the independent jet variables in
$(\mathcal{ E})$:
\begin{equation}\label{eee214}
U_\alpha^j=F_\alpha^j\left(x,u,(U_{\beta(q)}^{j(q)})_{1\leq q \leq p}\right),
\end{equation}
for $(j, \alpha) \neq (j(1), \beta(1)), \dots, (j(p), \beta (p))$ and
$j=1, \dots, m$, $\vert \alpha \vert \leq \kappa$. Clearly, the
natural coordinates on the submanifold $\Delta_\mathcal{E}$ of
$\mathcal{J}_{n,m}^\kappa$ are the $n+m+p$ coordinates
\begin{equation}\label{eee215}
\left(x,u,(U_{\beta(q)}^{j(q)})_{1\leq q \leq p}\right).
\end{equation}
Let $h=h(x,u)$ be a local $\K$-analytic diffeomorphism of $\K^{n+m}$
 close to the identity mapping and let
$\pi_\kappa: \mathcal{J}_{ n,m}^\kappa\to \K^{n+m}$ be the canonical
projection. According to~\cite{ol1986} (Chapter~2) there exists a
unique lift $h^{(\kappa)}$ of $h$ to $\mathcal{J}_{ n,m}^\kappa$ such that
$\pi_\kappa \circ h^{(\kappa)}= h\circ \pi_\kappa$. The components of
$h^{(\kappa)}$ may be computed by means of universal
combinatorial formulas and they are rational functions of the jet
variables~(\ref{eee212}), their coefficients being partial
derivatives of the components of $h$, {\it see} for instance \S3.3.5
of~\cite{bk1989}. By definition,
$h$ is a local symmetry of $(\mathcal{E})$ if $h$ transforms the
graph of every local solution of $(\mathcal{E})$ into the graph of
another local solution of $(\mathcal{E})$. This definition seems to be
rather uneasy to handle, because of the abstract quantification of
``every local solution'', but we have the following concrete
characterization for $h$ to be
a local symmetry of $(\mathcal{E})$, {\it cf.} Chapter~2 in~\cite{ol1986}.

\begin{lemma}\label{llem21}
The following conditions are equivalent:
\begin{itemize}
\item[{\bf (1)}]
The local transformation $h$ is a local symmetry of $(\mathcal{E})$.
\item[{\bf (2)}]
Its $\kappa$-th prolongation $h^{(\kappa)}$ is a local
self-transformation of the skeleton $\Delta_\mathcal{E}$ of
$(\mathcal{E})$.
\end{itemize}
\end{lemma}

These considerations have an infinitesimal version. Indeed, let
$X=\sum_{l=1}^n\, Q^l(x,u)\, \partial/\partial x_l+\sum_{j=1}^m\,
R^j(x,u)\, \partial/\partial u^j$ be a local vector field with
$\K$-analytic coefficients which is defined in a neighbourhood of the
origin in $\K^{n+m}$. Let $s\in\K$ and consider the flow of $L$ as the
one-parameter family $h_s(x,u):=\exp(s \, X)(x,u)$ of local
transformations.  We recall that $X$ is an {\sl infinitesimal
symmetry} of $(\mathcal{E})$ if for every small $s\in\K$, the mapping
$h_s(x,u):=\exp(s\, X)(x,u)$ is a local symmetry of
$(\mathcal{E})$. By differentiating with respect to $s$ the
$\kappa$-th prolongation $(h_s)^{(\kappa)}$ of $h_s$ at $s=0$, we
obtain a unique vector field $X^{(\kappa)}$ on the $\kappa$-th jet
space, called the {\sl $\kappa$-th prolongation of $X$} and which
satisfies $(\pi_k)_*(X^{(\kappa)})=X$. In Subsections~3.1 and~3.2
below, we shall analyze the combinatorial formulas for the
coefficients of $X^{(\kappa)}$, since they will be needed to prove
Theorem~\ref{thm1}.

Let $X_\mathcal{E}$ be the projection to the restricted jet space
$\K^{m+n+p}$, equipped with the coordinates~(\ref{eee215}), of the
restriction of $X^{(\kappa)}$ to $\Delta_{\mathcal{E}}$, namely
\begin{equation}\label{eee216}
X_\mathcal{E}:=(\pi_{\kappa,p})_*(X^{(\kappa)}\vert_{\Delta_{\mathcal{E}}}).
\end{equation}
The following Lemma, called the {\sl Lie criterion}, is the concrete
characterization for $X$ to be an infinitesimal symmetry of
$(\mathcal{E})$ and is a direct corollary of Lemma~\ref{llem21}, {\it
cf.} Chapter~2 in~\cite{ol1986}. This criterion will be central in the
next Sections~3, 4 and 5.

\begin{lemma}\label{llem22}
The following conditions are equivalent:
\begin{itemize}
\item[{\bf (1)}]
The vector field $X$ is an infinitesimal symmetry of $(\mathcal{E})$.
\item[{\bf (2)}]
Its $\kappa$-th prolongation $X^{(\kappa)}$ is tangent to
the skeleton $\Delta_\mathcal{E}$.
\end{itemize}
\end{lemma}

We denote by $\mathfrak{Sym}(\mathcal{E})$ the set of infinitesimal
symmetries of ($\mathcal E$). Since it may be easily checked that
$(cX+dY)^{(\kappa)}=cX^{(\kappa)}+ dY^{(\kappa)}$ and that
$[X^{(\kappa)}, Y^{(\kappa)}]= \left( [X,Y] \right)^{(\kappa)}$, {\it
see} Theorem~2.39 in~\cite{ol1986}, it follows from Lemma~\ref{llem22}
{\bf (2)} that $\mathfrak{Sym}(\mathcal{E})$ is a Lie algebra of
locally defined vector fields. Our main question in this section is
the following: {\it under which natural conditions is
$\mathfrak{Sym}(\mathcal{E})$ finite-dimensional~?}

\begin{example}\label{ex22}
{\rm 
We observe that the Lie algebra $\mathfrak{Sym}(\mathcal{E})$ of the
system ($\mathcal E$) presented in Example~\ref{ex21} is
infinite-dimensional, since it includes all vector fields of the form
$X=Q^2(x_1,x_2,u)\, \partial /\partial x_2$, as may be verified. As we
will argue in Proposition~\ref{pprop22} below, this phenomenon is
typical, the main reason lying in the first order relation
$u_{x_2}=0$.
}
\end{example}

By analyzing the construction of the submanifold of solutions
$\mathcal{M}$ associated to the system $(\mathcal{E})$, we may
establish the following correspondence (we shall not develop its
proof).

\begin{proposition}\label{pprop21}
To every infinitesimal symmetry $X=\sum_{l=1}^n\, 
Q^l(x,u)\, \partial/\partial x_l+\sum_{j=1}^m\,
R^j(x,u)\, \partial /\partial u^j$ of $(\mathcal{E})$, 
there corresponds a unique vector field of the form
\begin{equation}\label{eee217}
\mathcal{X}=
\sum_{j=1}^m\, 
\Pi^j(\nu,\chi)\,
\frac{\partial }{\partial \nu^j}+
\sum_{q=1}^p\, 
\Lambda^q(\nu,\chi)\, 
\frac{\partial }{\partial \chi_q},
\end{equation}
whose coefficients depend only on the parameters $(\nu,\chi)$, such
that $X+\mathcal{X}$ is tangent to the submanifold of solutions
$\mathcal{M}$.
\end{proposition}

This leads us to define the Lie algebra $\mathfrak{Sym}(\mathcal{M})$
of vector fields of the form
\begin{equation}\label{eee218}
\sum_{l=1}^n\, 
Q^l(x,u)\, \frac{\partial}{\partial x_l}+\sum_{j=1}^m\,
R^j(x,u)\, \frac{\partial}{\partial u^j}+
\sum_{j=1}^m\, 
\Pi^j(\nu,\chi)\,
\frac{\partial }{\partial \nu^j}+
\sum_{q=1}^p\, 
\Lambda^q(\nu,\chi)\, 
\frac{\partial }{\partial \chi_q}
\end{equation}
which are tangent to $\mathcal{M}$. We shall say that the submanifold
$\mathcal{M}$ is {\sl degenerate} if there exists a nonzero vector
field of the form $X=\sum_{l=1}^n\, Q^l(x,u)\, \partial/\partial
x_l+\sum_{j=1}^m\, R^j(x,u)\, \partial /\partial u^j$ which is tangent
to $\mathcal{M}$, which means that the corresponding $\mathcal{X}$
part is zero. In this case, we claim that
$\mathfrak{Sym}(\mathcal{M})$ is infinite dimensional. Indeed there
exists then a nonzero vector field $T~=~\sum_{l=1}^n Q^l(x,u) \partial
/\partial x_l\ + \sum_{j=1}^m R^j(x,u)\partial / \partial u^j$ tangent
to $\mathcal M$. Consequently, for every $\K$-analytic function
$A(x,u)$, the vector field $A(x,u)\, T$ belongs to
$\mathfrak{Sym}(\mathcal{M})$, hence $\mathfrak{Sym}(\mathcal{M})$ is
infinite dimensional.
 
By developing the dual defining functions of
$\mathcal{M}$ with respect to the powers of $\chi$, we may write
\begin{equation}\label{eeee20}
\nu^j=\Omega_j^*(\chi,x,u)=
\sum_{\gamma\in\N^p}\, \chi^\gamma \Omega_{j,\gamma}^*(x,u),
\end{equation}
where the functions $\Omega_{j,\gamma}^*(x,u)$ are
$\K$-analytic in a neighbourhood of the origin,
we may formulate a criterion for $\mathcal{M}$ to be non degenerate with
respect to the variables (whose proof is skipped).

\begin{proposition}\label{pprop22}
The submanifold $\mathcal{M}$ is {\rm not} degenerate with respect to
the variables if and only if there exists an integer $k$ such that the
generic rank of the local $\K$-analytic mapping
\begin{equation}\label{eeee21}
(x,u)\longmapsto 
\left(
\Omega_{j,\gamma}^*(x,u)
\right)_{1\leq j\leq m, \ \gamma\in\N^p, \ \vert \gamma \vert \leq k} 
\end{equation}
is equal to $n+m$.
\end{proposition}

Seeking for conditions which insure that
$\mathfrak{Sym}(\mathcal{M})$ is finite-dimensional, it is therefore
natural to assume that the generic rank of the mapping~(\ref{eeee21})
is equal to $n+m$. Furthermore, to simplify the presentation, we shall
assume that the {\it rank at $(x,u)=(0,0)$} (not only the generic
rank) {\it of the mapping~(\ref{eeee21}) is equal to $n+m$ for $k$
large enough}. This is a ``Zariski-generic'' assumption. Coming back
to~(\ref{eee211}), we observe that this means exactly that
$\mathcal{M}$ is solvable with respect to the variables.  Then we
denote by $\ell_0^*$ the smallest integer $k$ such that the rank at
$(x,u)=(0,0)$ of the mapping~(\ref{eeee21}) is equal to $n+m$ and we
say that $\mathcal{M}$ is {\sl $\ell_0^*$-solvable with respect to the
variables}.  Also, we denote by $\ell_0$ the integer $\max_{1\leq
q\leq p}\, \vert \beta(q) \vert$ and we say that $\mathcal{M}$ is {\sl
$\ell_0$-solvable with respect to the parameters}.

\subsection{Fundamental isomorphism between 
$\mathfrak{Sym}(\mathcal{E})$ and $\mathfrak{Sym}(\mathcal{M})$} In
the remainder of this Section~2, we shall assume that $\mathcal{M}$ is
$\ell_0$-solvable with respect to the parameters and
$\ell_0^*$-solvable with respect to the variabes. In this case,
viewing the variables $(\nu^1,\dots,\nu^m)$ in the dual equations
$\nu^j=\Omega_j^*(\chi,x,u)$ of $\mathcal{M}$ as a mapping of $\chi$
with (dual) ``parameters'' $(x,u)$ and proceeding as in
Subsection~2.2, we may construct a {\sl dual system of completely
integrable partial differential equations} of the form
\vskip 0,1cm
\noindent {$(\mathcal{E}^*)$}\centerline{
$\nu_{\chi^\gamma}^j(\chi)=
G_\gamma^j\left(
\chi,\nu(\chi),(\nu_{\chi^{\delta(l)}}^{j(l)}(\chi))_{
1\leq l\leq n}
\right),$
}

\vskip 0,1cm
\noindent where $(j,\gamma)\neq (j(1),\delta(1)),\dots,(j(n),\delta(n))$.
This system has its own infinitesimal symmetry Lie algebra
$\mathfrak{Sym}(\mathcal{E}^*)$. 

\begin{theorem}\label{tthm22}
If $\mathcal{M}$ is both solvable with respect to the
parameters and solvable with respect to the variables, 
we have the following two isomorphisms:
\begin{equation}\label{eee222}
\mathfrak{Sym}(\mathcal{E}) \cong \mathfrak{Sym}(\mathcal{M})\cong
\mathfrak{Sym}(\mathcal{E}^*),
\end{equation}
namely $X \longleftrightarrow X+\mathcal{X} \longleftrightarrow \mathcal{X}$.
\end{theorem}

In Subsection~2.10 below, we shall introduce a second geometric condition 
which is in general necessary for $\mathfrak{Sym}(\mathcal{M})$ 
to be finite-dimensional.

\subsection{Local (pseudo)group ${\rm Sym}(\mathcal{M})$ of point 
transformations of $\mathcal{M}$}
We shall study the geometry of a local $\K$-analytic submanifold
$\mathcal{ M}$ of $\K^{ n+2m+p}$ whose equations and dual equations
are of the form
\begin{equation}\label{eee223}
\left\{
\aligned
u^j=
& \
\Omega_j(x,\nu,\chi), \ \ \ \ \ \ \ \ \ j=1,\dots,m,\\
\nu^j=
& \
\Omega_j^*(\chi,x,u), \ \ \ \ \ \ \ \ \ j=1,\dots,m.
\endaligned\right.
\end{equation}
Let $t:=(x,u)\in\K^{n+m}$ and $\tau:=(\nu,\chi)\in\K^{n+m}$. We are
interested in describing the set of local $\K$-analytic
transformations of the space $\K^{n+2m+p}$ which are of the specific
form
\begin{equation}\label{eee224}
(t,\tau)\longmapsto (h(t),\phi(\tau)),
\end{equation}
and which stabilize $\mathcal{ M}$, in a neighborhood of the
origin. We denote the local Lie pseudogroup of such transformations
(possibly infinite-dimensional) by ${\rm Sym}( \mathcal{
M})$. Importantly, each transformation of ${\rm Sym} (\mathcal{ M})$
stabilize both the sets $\{t= ct.\}$ and the sets $\{\tau= ct.\}$. Of
course, the Lie algebra of ${\rm Sym}(\mathcal{ M})$ coincides with
$\mathfrak{ Sym} ( \mathcal{ M})$ defined above.

\subsection{Fundamental pair of foliations on $\mathcal{M}$}
Let $p_0\in\K^{n+2m+p}$ be a fixed point of coordinates
$(t_{p_0},\tau_{p_0})$.  Firstly, we observe that the intersection
$\mathcal{M}\cap \{\tau=\tau_{p_0}\}$ consists of the $n$-dimensional
$\K$-analytic submanifold of equation $u=\Omega(x,\tau_{p_0})$. As
$\tau_{p_0}$ varies, we obtain a local $\K$-analytic foliation of
$\mathcal{M}$ by $n$-dimensional submanifolds. Let us denote this
first foliation by $\mathcal{F}_p$ and call it the {\sl foliation of
$\mathcal{M}$ with respect to parameters}.  Secondly, and dually, we
observe that the intersection $\mathcal{M}\cap \{t=t_{p_0}\}$ consists
of the $p$-dimensional $\K$-analytic submanifold of equation
$\nu=\Omega^*( \chi,t_{p_0})$. As $t_{p_0}$ varies, we obtain a local
$\K$-analytic foliation of $\mathcal{M}$ by $p$-dimensional
submanifolds. Let us denote this second foliation by $\mathcal{F}_v$
and call it the {\sl foliation of $\mathcal{M}$ with respect to the
variables}. We call $(\mathcal{F}_p, \mathcal{F}_v)$ the 
{\sl fundamental pair of foliations on $\mathcal{M}$}.

\subsection{Covering property of the fundamental pair of foliations}
We wish to formulate a geometric condition which says that starting
from the origin in $\mathcal{M}$ and following alternately the leaves
of $\mathcal{F}_p$ and the leaves of $\mathcal{F}_v$, we cover a
neighborhood of the origin in $\mathcal{M}$. Let us introduce two
collections $(\mathcal{L}_k)_{1\leq k\leq n}$ and
$(\mathcal{L}_q^*)_{1\leq q\leq p}$ of vector fields whose integral
manifolds coincide with the leaves of $\mathcal{F}_p$ and
$\mathcal{F}_v$:
\begin{equation}\label{eee225}
\left\{
\aligned
\mathcal{L}_k:=
\frac{\partial }{\partial x_k}+
\sum_{j=1}^m \, 
\frac{\partial \Omega_j}{\partial x_k} (x,\nu,\chi)\, 
\frac{\partial }{\partial u^j}, 
\ \ \ \ \ \ k=1,\dots,n, \\
\mathcal{L}_q^*:=
\frac{\partial }{\partial \chi_q}+
\sum_{j=1}^m\, 
\frac{\partial \Omega_j^*}{\partial \chi_q}(\chi,x,u)\,
\frac{\partial }{\partial \nu^j},
\ \ \ \ \ \ k=1,\dots,n.
\endaligned\right. 
\end{equation}
Let $p_0$ be a fixed point in $\mathcal{M}$ of coordinates
$(x_{p_0},u_{p_0},\nu_{p_0},\chi_{p_0})\in\K^{n+2m+p}$, let
$x_1:=(x_{1,1},\dots,x_{1,n})\in\K^n$ be a ``multitime'' parameter and
define the multiple flow map
\begin{equation}\label{eee226}
\left\{
\aligned
\mathcal{L}_{x_1}(x_{p_0},u_{p_0},\nu_{p_0},\chi_{p_0}):=
& \
\exp(x_1\mathcal{L})(p_0):=
\exp(x_{1,n}\mathcal{L}_n(\cdots(\exp(x_{1,1}\mathcal{L}_1(p_0)))\cdots)):=\\
:=
& \
(x_{p_0}+x_1,\Omega(x_{p_0}+x_1,\nu_{p_0},\chi_{p_0}),\nu_{p_0},\chi_{p_0}).
\endaligned\right.
\end{equation}
Similarly, for $\chi=(\chi_{1,1},\dots,\chi_{1,p})\in\K^p$, define the
multiple flow map
\begin{equation}\label{eee227}
\mathcal{L}_{\chi_1}^*(x_{p_0},u_{p_0},\nu_{p_0},\chi_{p_0}):=
(x_{p_0},u_{p_0},\Omega^*(\chi_{p_0}+\chi_1,
x_{p_0},u_{p_0}),\chi_{p_0}+\chi_1).
\end{equation}
We may define now the mappings which
correspond to start from the origin and to move alternately along the
two foliations $\mathcal{F}_p$ and $\mathcal{F}_v$. If the first
movement consists in moving along the foliation $\mathcal{F}_v$, we
define 
\begin{equation}\label{eee228}
\left\{
\aligned
\Gamma_1(x_1):=
& \
\mathcal{L}_{x_1}(0), \\
\Gamma_1(x_1,\chi_1):=
& \
\mathcal{L}_{\chi_1}^*(\mathcal{L}_{x_1}(0)), \\
\Gamma_3(x_1,\chi_1,x_2):=
& \
\mathcal{L}_{x_2} (\mathcal{L}_{\chi_1}^*(\mathcal{L}_{x_1}(0))), \\
\Gamma_4(x_1,\chi_1,x_2,\chi_2):=
& \
\mathcal{L}_{\chi_2}^*(
\mathcal{L}_{x_2} (\mathcal{L}_{\chi_1}^*(\mathcal{L}_{x_1}(0)))).
\endaligned\right.
\end{equation}
Generally, we may define the maps $\Gamma_k ([x \chi]_k)$, where $[x
\chi]_k= (x_1, \chi_1, x_2, \chi_2, \dots)$ with exactly $k$ terms and
where each $x_l$ belongs to $\K^n$ and each $\chi_l$ belongs to
$\K^p$. On the other hand, if the first movement consists in moving
along the foliation $\mathcal{ F}_p$, we start with $\Gamma_1^*(
\chi_1):= \mathcal{ L}_{ \chi_1}^* (0)$, $\Gamma_2^*( \chi_1, x_1):=
\mathcal{ L}_{x_1}( \mathcal{ L}_{ \chi_1}^* (0))$, {\it etc.}, and
generally we may define the maps $\Gamma_k^* ([\chi x]_k)$, where
$[\chi x]_k=( \chi_1, x_1, \chi_2, x_2, \dots)$, with exactly $k$
terms.  The range of both maps $\Gamma_k$ and $\Gamma_k^*$ is
contained in $\mathcal{M}$. We call $\Gamma_k$ the {\sl $k$-th chain}
and $\Gamma_k^*$ the {\sl $k$-th dual chain}.

\begin{definition}\label{ddef21}
{\rm 
The pair of foliations $(\mathcal{F}_p,\mathcal{F}_v)$ is called {\sl
covering at the origin} if there exists an integer $k$ such that the
generic rank of $\Gamma_k$ is (maximal possible) equal to $\dim_\K \,
\mathcal{M}$.  Since the dual $(k+1)$-th chain $\Gamma_{k+1}^*$ for
$\chi_1=0$ identifies with the $k$-th chain $\Gamma_k$, it follows
that the same property holds for the dual chains.
}
\end{definition}

In terms of Sussmann's approach~\cite{sus1973}, this means that
the {\sl local orbit} of the two systems of vector fields
$(\mathcal{L}_k)_{1\leq k\leq n}$ and $(\mathcal{L}_q^*)_{1\leq q \leq
p}$ is of maximal dimension.  Reasoning as in~\cite{sus1973} (using
the so-called {\sl backward trick} in Control Theory, {\it see}
also~\cite{me2003}), it may be shown that there exists the
smallest {\it even}\, integer $2\mu_0$ such that the ranks of the two
maps $\Gamma_{2\mu_0}$ and $\Gamma_{2\mu_0}^*$ at the origin (not only
their generic rank) in $\K^{n\mu_0+p\mu_0}$ are both equal to $\dim_\K
\, \mathcal{M}$. This means that $\Gamma_{2\mu_0}$ and
$\Gamma_{2\mu_0}^*$ are submersive onto a neighborhood of the origin
in $\mathcal{M}$. We call $\mu_0$ the {\sl type of the pair of
foliations $(\mathcal{F}_p, \mathcal{F}_v)$}. It may also be established
that $\mu_0\leq m+2$.

\def\theexample{2.46}\begin{example}\label{ex23}
{\rm
We give an example of a submanifold which is both $1$-solvable with
respect to the parameters and with respect to the variables but whose pair of
foliations is not covering: with $n=1$, $m=2$ and $p=1$, this is given
by the two equations $u^1=\nu^1$, $u^2=\nu^2+x\chi_1$. Then
$\mathfrak{Sym}(\mathcal{M})$ is infinite-dimensional since it
contains the vector fields $a(u^1)\, \partial/\partial u^1+
a(\nu^1)\, \partial /\partial \nu^1$, where $a$ is an arbitrary 
$\K$-analytic function. For this reason, we shall assume in the sequel that
the pair of foliations $(\mathcal{F}_p, \mathcal{F}_v)$ is
covering at the origin.
}
\end{example}

\subsection{Estimate on the dimension of the local symmetry group of the
submanifold of solutions} We may now formulate the main theorem of
this section, which shows that, under suitable nondegeneracy
conditions, ${\rm Sym}(\mathcal{M})$ is a finite dimensional local Lie
group of local transformations.  If $t\in\K^{n+m}$, we denote by
$\vert t\vert:=\max_{1\leq k\leq n+m}\, \vert t_k \vert$. If $(h,\phi)
\in Sym(\mathcal M)$ we denote by $J_t^kh(0)$ the $k$-th order jet of
$h$ at the origin and by $J_\tau^k\phi(0)$ the $k$-th order jet of
$\phi$ at the origin. Also, we shall assume that $\mathcal{M}$ is
either $\K$-algebraic or $\K$-analytic. Of course, the
$\K$-algebraicity of the submanifold of solutions does not follow from
the $\K$-algebraicity of the right hand sides $F_\alpha^j$ of the
system of partial differential equations $(\mathcal{E})$.

\begin{theorem}\label{tthm23}
Assume that the $\K$-algebraic 
or $\K$-analytic submanifold of solutions $\mathcal{M}$
of the completely integrable
system of partial differential 
equations $(\mathcal{E})$ is both $\ell_0$-sovable with respect to the
parameters and $\ell_0^*$-solvable with respect to the
variables. Assume that the fundamental pair of foliations
$(\mathcal{F}_p,\mathcal{F}_v)$ is covering at the origin and let
$\mu_0$ be its type at the origin. Then there exists $\varepsilon_0>0$
such that for every $\varepsilon$ with $0<\varepsilon <
\varepsilon_0$, the following four properties hold:
\begin{itemize}
\item[{\bf (a)}]
The (pseudo)group ${\rm Sym}(\mathcal{M})$ of local $\K$-analytic
diffeomorphisms defined for $\{(t,\tau)\in\K^{n+2m+p}: \, \vert t
\vert < \varepsilon, \ \vert \tau \vert < \varepsilon\}$ which are of
the form $(t,\tau)\mapsto (h(t),\phi(\tau))$ and which stabilize
$\mathcal{M}$ is a local Lie pseudogroup of transformations of
finite dimension $d\in\N$.
\item[{\bf (b)}]
Let $\kappa_0:=\mu_0(\ell_0+\ell_0^*)$. Then there exist two
$\K$-algebraic or $\K$-analytic mappings $H_{\kappa_0}$ and
$\Phi_{\kappa_0}$ which depend only on $\mathcal{M}$ and which may be
constructed algorithmically by means of the defining equations of
$\mathcal{M}$ such that every element $(h,\phi)\in{\rm
Sym}(\mathcal{M})$, sufficiently close to the identity mapping,
may be represented by
\begin{equation}\label{eee229}
\left\{
\aligned
h(t) =
& \
H_{\kappa_0}(t,\, J_t^{\kappa_0} h(0)), \\
\phi(\tau)=
& \
\Phi_{\kappa_0}(\tau,\, J_\tau^{\kappa_0} \phi(0)).
\endaligned\right.
\end{equation}
Consequently, every element of\, ${\rm Sym}(\mathcal{M})$ is
uniquely determined by its $\kappa_0$-th jet at the origin
and the dimension $d$ of the Lie algebra $\mathfrak{Sym}(\mathcal{M})$
is bounded by the number of components of the vector
$(J_t^{\kappa_0}h(0), \, J_\tau^{\kappa_0} \phi(0))$, namely 
we have
\begin{equation}\label{eee230}
\dim_\K\, \mathfrak{Sym}(\mathcal{E})=
\dim_\K \, \mathfrak{Sym}(\mathcal{M})\leq
(n+m)\, C_{n+m+\kappa_0}^{\kappa_0}+
(m+p)\, C_{m+p+\kappa_0}^{\kappa_0}.
\end{equation}
\item[{\bf (c)}]
There exists $\varepsilon'$ with $0< \varepsilon' < \varepsilon$ and a
$\K$-algebraic or $\K$-analytic mapping $( H_\mathcal{ M}, \,
\Phi_\mathcal{ M})$ which may be constructed algorithmically by means
of the defining equations of $\mathcal{M}$, defined in a
neighbourhood of the origin in $\K^{ n+2m+p} \times \K^d$ with values
in $\K^{ n+2m+p}$ and which satifies $(H_\mathcal{ M} (t,0), \,
\Phi_\mathcal{ M} (\tau,0)) \equiv (t,\tau)$, such that every element
$(h, \phi) \in {\rm Sym} (\mathcal{ M})$ defined on the set
$\{(t,\tau)\in\K^{ n+2m+p}: \, \vert t \vert < \varepsilon', \ \vert
\tau \vert < \varepsilon' \}$, sufficiently close to the identity
mapping and stabilizing $\mathcal{ M}$ may be represented as $(h(t),
\phi (\tau)) \equiv (H_\mathcal{M} (t,s_{h,\phi}), \, \Phi_\mathcal{M}
(\tau, s_{h,\phi}))$ for a unique element $s_{h,\phi} \in \K^d$
depending on the mapping $(h,\phi)$.
\item[{\bf (d)}]
The mapping $(t,\tau,s)\longmapsto (H_\mathcal{M}(t,s),
\Phi_\mathcal{M}(\tau,s))$ defines a local $\K$-algebraic or
$\K$-analytic Lie group of local $\K$-algebraic or $\K$-analytic
transformations stabilizing $\mathcal{M}$.
\end{itemize}
\end{theorem}

\subsection{Applications}
The proof of Theorem~\ref{tthm23}, which possesses strong similarities
with the proof of Theorem~4.1 in~\cite{gm2002a}, will not be
presented. It seems that Theorem~\ref{tthm23}, together with the
argumentation on the necessity of assumptions that $\mathcal{M}$ be
solvable with respect to the variables and that its fundamental pair
of foliations be covering, is a new result about the
finite-dimensionality of a completely integrable system of partial
differential equations having an arbitrary number of independent and
dependent variables. The main interest lies in the fact that we obtain
the algorithmically constructible representation
formula~(\ref{eee229}) together with the local Lie group structure
mapping $(H_{\mathcal{M}}, \, \Phi_{\mathcal{M}})$. In particular, we
get as a corollary that every transformation $(h(t),\, \phi(\tau))$
given by a formal power series (not necessarily convergent) is as
smooth as the applications $(H_{\kappa_0}, \, \Phi_{\kappa_0})$ are,
namely every formal element of ${\rm Sym}(\mathcal{M})$ is necessarily
$\K$-algebraic or $\K$-analytic. As a counterpart of its generality,
Theorem~3 does not provide optimal bounds, as shows the following
illustration.

\begin{example}\label{ex24}
{\rm 
Let $n=m=1$, let $\kappa\geq 3$ and let $(\mathcal{ E})$ denote the
ordinary differential equation $u_{ x^\kappa } (x) = F(x, u(x),
u_x(x), \dots, u_{x^{ \kappa-1}} (x))$. Then the submanifold of
solutions $\mathcal{ M}$ is of the form $u= \nu+ x\chi_1+ \cdots+ x^{
\kappa-1} \chi_{\kappa-1 }+ {\rm O}( \vert x \vert^\kappa) + {\rm O}
(\vert \chi \vert^2)$.  It may be checked that $\ell_0= \kappa-1$,
$\ell_0^* =1$ and $\mu_0=3$, hence $\kappa_0=3\kappa$. Then the
dimension estimate in~(\ref{eee230}) is: $\dim_\K \,
\mathfrak{Sym}(\mathcal{E})\leq 2\, C_{2+3\kappa}^{3\kappa}+ \kappa\,
C_{4\kappa}^{3\kappa}$. This bound is much larger than the optimal
bound $\dim_\K \, \mathfrak{Sym}(\mathcal{E}) \leq \kappa + 4$ due to
S.~Lie ({\it cf.}~\cite{enlie1890}; {\it see} also the case $n=m=1$ of
Theorem~\ref{thm1}).
}
\end{example}

Untill now we focused on providing the set of Lie symmetries of a
general system of partial differential equations with a local Lie
group structure. As a byproduct we obtained the (non optimal)
dimensional upper bound (\ref{eee230}) of Theorem~\ref{tthm23}. In the
next Sections~3, 4 and 5, using the classical Lie algorithm based on
the Lie criterion (see Lemma~\ref{llem22}), we provide an optimal
bound for some specific systems of partial differential equations,
answering an open problem raised in~\cite{ol1995} page~206.

\section{Lie theory for partial differential equations}

\subsection{Prolongation of vector fields to the jet spaces}
Consider the following $\K$-analytic system $(\mathcal E)$ of non
linear partial differential equations:
\begin{equation}\label{e21}
u_{x_{k_1}\cdots x_{k_\kappa}}^j(x)= F_{k_1,\dots,k_\kappa}^j\left(x,
u(x), u_{x_{l_1}}^{i_1}(x),\dots,u_{x_{l_1}\cdots
x_{l_{\kappa-1}}}^{i_1}(x)\right),
\end{equation}
where $1\leq k_1\leq \cdots \leq k_\kappa\leq n$, $1 \leq j \leq m$,
and $F_{k_1,\dots,k_\kappa}^j$ are analytic functions of $n+m\,
C_{n+\kappa-1}^{\kappa-1}$ variables, defined in a neighbourhood of
the origin. We assume that $(\mathcal E)$ is {\sl completely
integrable}. The Lie theory consists in studying the {\sl
infinitesimal symmetries} $X=\sum_{l=1}^n \, Q^l(x,u)\,
\partial/\partial x_l+ \sum_{j=1}^m\, R^j(x,u)\, \partial /\partial
u^j$ of $(\mathcal{E})$. Consider the {\sl skeleton} of
$(\mathcal{E})$, namely the complex subvariety $\Delta_\mathcal{E}$ of
codimension $m \, C_{\kappa+n-1}^\kappa$ in the jet space
$\mathcal{J}_{n,m}^\kappa$, defined by
\begin{equation}\label{e23}
U_{k_1,\dots,k_\kappa}^j=F_{k_1,\dots,k_\kappa}^j\left(x,
u,U_{l_1}^{i_1},\dots,
U_{l_1,\dots,l_{\kappa-1}}^{i_1}\right),
\end{equation}
where $j,i_1=1,\dots,m$ and $k_1,\dots,k_\kappa,l_1,\dots,
l_{\kappa-1}=1,\dots,n$. For $k=1,\dots,n$ let $D_{k}$ be the $k$-th
{\sl operator of total differentiation}, characterized by the property
that for every integer $\lambda \geq 2$ and for every analytic
function $P=P(x,u,U_{l_1}^{i_1},\dots,
U_{l_1,\dots,l_{\lambda-1}}^{i_1})$ defined in the jet space
$\mathcal{J}_{n,m}^{\lambda-1}$, the operator $D_{k}$ is the unique
formal infinite differential operator satisfying the relation
\begin{equation}\label{e24}
\left\{
\aligned
{}
&
[D_{k}P]\left(x,u(x),
u_{x_{l_1}}^{i_1}(x),\dots,
u_{x_{l_1}\cdots x_{l_{\lambda-1}}}^{i_1}(x)
\right)\equiv \\
& \ \ \ \ \ \ \ \ \ \ \ \ \ \ \ \ \ \ 
\frac{\partial }{\partial x_{k}}\left[
P\left(x,u(x),
u_{x_{l_1}}^{i_1}(x),\dots,
u_{x_{l_1}\cdots x_{l_{\lambda-1}}}^{i_1}(x)
\right)
\right].
\endaligned\right.
\end{equation}
Note that this identity involves only the troncature of $D_k$ to order
$\lambda$, denoted by $D_{k}^\lambda$, and defined by
\begin{equation}\label{e25}
\left\{
\aligned
{}
&
D_{k}^\lambda:=
\frac{\partial }{\partial x_{k}}+
\sum_{i_1=1}^m\, U_{k}^{i_1} \frac{\partial}{\partial u^{i_1}}+
\sum_{i_1=1}^m\, \sum_{l_1=1}^n\, 
U_{k,l_1}^{i_1}\, 
\frac{\partial }{\partial U_{l_1}^{i_1}}+\cdots+\\
& \ \ \ \ \ \ \ \ \ \ \ \ \ \ \ \ \ \ \ \ 
+
\sum_{i_1=1}^m\, \sum_{l_1,\dots,l_{\lambda-1}=1}^n\,
U_{k,l_1,\dots,l_{\lambda-1}}^{i_1}\, 
\frac{\partial }{\partial U_{l_1,\dots,l_{\lambda-1}}^{i_1}}.
\endaligned\right.
\end{equation}
According to Theorem~2.36 of~\cite{ol1986}, the {\sl prolongation of
order $\kappa$ of a vector field $X= \sum_{l=1}^n\, Q^l(x,u)\,
\partial/\partial x_l+ \sum_{j=1}^m\, R^j(x,u)\, \partial /\partial
u^j$}, denoted by $X^{(\kappa)}$, is the unique vector field on the
space $\mathcal{J}_{n,m}^\kappa$ of the form
\begin{equation}\label{e26}
\left\{
\aligned
{}
&
X^{(\kappa)}=
X+
\sum_{j=1}^m\, 
\sum_{k_1=1}^n\,
{\bf R}_{k_1}^j\, 
\frac{\partial }{\partial U_{k_1}^j}+
\sum_{j=1}^m\, \sum_{k_1,k_2=1}^n\, 
{\bf R}_{k_1,k_2}^j\, 
\frac{\partial }{\partial U_{k_1,k_2}^j}+\cdots+\\
& \ \ \ \ \ \ \ \ \ \ \ \ \ \ \ \ \ \ \ \ \ \ \ \
+\sum_{j=1}^m\, 
\sum_{k_1,\dots,k_\kappa=1}^n\, 
{\bf R}_{k_1,\dots,k_\kappa}^j\, 
\frac{\partial }{\partial U_{k_1,k_2,\dots,k_\kappa}^j},
\endaligned\right.
\end{equation}
corresponding to the infinitesimal action of the flow of $X$ on the
jets of order $\kappa$ of the graphs of maps $u=u(x)$, and whose
coefficients are computed recursively by the formulas
\begin{equation}\label{e27}
\left\{
\aligned
{\bf R}_{k_1}^j:= 
& \
D_{k_1}^1(R^j)-\sum_{l_1=1}^n\, 
D_{k_1}^1(Q^{l_1})\, U_{l_1}^j,\\
{\bf R}_{k_1,k_2}^j:= 
& \
D_{k_2}^2 ({\bf R}_{k_1}^j)-
\sum_{l_1=1}^n\, 
D_{k_2}^1(Q^{l_1})\, 
U_{k_1,l_1}^j,\\
\cdots \cdots
& \cdots \cdots \cdots \cdots \cdots \cdots \cdots\\
{\bf R}_{k_1,k_2,\dots,k_\lambda}^j:=
& \
D_{k_\lambda}^\lambda({\bf R}_{k_1,\dots,k_{\lambda-1}}^j)-
\sum_{l_1=1}^n\, D_{k_\lambda}^1(Q^{l_1})\, 
U_{k_1,\dots,k_{\lambda-1},l_1}^j.
\endaligned\right.
\end{equation}
For a better comprehension of the general computation, let us start by
computing $R^\kappa$ in the case $n=m=1$. 

\subsection{Computation of $R^\kappa$ when $n=m=1$} A direct
application of the preceding formulas leads to the following classical
expressions:
\begin{equation}\label{e28}
\left\{
\aligned
{\bf R}^1= & \ R_x+
[R_u-Q_x]\,U^1+
[-Q_u] \, (U_1)^2. \\
{\bf R}^2 = & \ R_{x^2}+
[2R_{xu}-Q_{x^2}]\,U^1+
[R_{u^2} -2Q_{xu}]\,(U^1)^2+
[-Q_{u^2}]\,(U^1)^3+\\
& \ \ \ \ \ \ \ +
[R_u-2Q_x] \, U^2+
[-3Q_u] \, U^1 U^2.\\
\endaligned\right.
\end{equation}
Observe that these expressions are polynomial in the jet variables,
their coefficients being differential expressions involving a partial
derivative of $R$ (with a positive integer coefficient) and a partial
derivative of $Q$ (with a negative integer coefficient). We have also:
\begin{equation}\label{e29}
\left\{
\aligned
{\bf R}^3= & \ R_{x^3}+
[3R_{x^2u}-Q_{x^3}]\,U^1+
[3R_{xu^2}-3Q_{x^2u}]\,(U^1)^2+\\
& \ \ \ \ \ \ \ +
[R_{u^3}-3Q_{xu^2}]\,(U^1)^3+
[-Q_{u^3}] \, (U^1)^4+
[3R_{xu}-3Q_{x^2}]\,U^2+\\
& \ \ \ \ \ \ \ +
[3R_{u^2} -9Q_{xu}]\,U^1U^2+
[-6Q_{u^2}]\,(U^1)^2 U^2+
[-3Q_u]\,(U^2)^2+\\
& \ \ \ \ \ \ \ +
[R_u -3Q_x]\,U^3+
[-4Q_u]\,U^1U^3.\\
{\bf R}^4= & \ R_{x^4}+
[4R_{x^3u}-Q_{x^4}]\,U^1+
[6R_{x^2u^2}-4Q_{x^3u}]\,(U^1)^2+\\
& \ \ \ \ \ \ \ +
[4R_{xu^3}-6Q_{x^2u^2}]\,(U^1)^3+
[R_{u^4}-4Q_{xu^3}]\,(U^1)^4+
[-Q_{u^4}]\,(U^1)^5+\\
& \ \ \ \ \ \ \ +
[6R_{x^2u}-4Q_{x^3}]\,U^2+
[12R_{xu^2}-18 Q_{x^2u}]\,U^1U^2+\\
& \ \ \ \ \ \ \ +
[6R_{u^3}-24Q_{xu^2}]\,(U^1)^2U^2+
[-10Q_{u^3}]\,(U^1)^3 U^2+\\
& \ \ \ \ \ \ \ +
[3R_{u^2}-12Q_{xu}]\,(U^2)^2+
[-15Q_{u^2}]\,U^1(U^2)^2+
[4R_{xu}-6Q_{x^2}]\,U^3+\\
& \ \ \ \ \ \ \ +
[4R_{u^2}-16 Q_{xu}]\,U^1U^3+
[-10Q_{u^2}]\,(U^1)^2U^3+
[-10Q_u]\,U^2U^3+\\
& \ \ \ \ \ \ \ +
[R_u-4Q_x]\,U^4+
[-5Q_u]\,U^1U^4.\\
\endaligned\right.
\end{equation}
Remark that all the brackets involved in equations~(\ref{e29}) are of
the form $[\lambda \, R_{x^au^{b+1}}-\mu \, Q_{x^{a+1} u^b}]$, where
$\lambda, \mu\in\N$ and $a, b\in\N$.

In what follows we will not need the complete form of $R^\kappa$ but
only the following partial form:
\begin{lemma}\label{lem21}
For $\kappa\geq 4$:
\begin{equation}\label{e210}
\left\{
\aligned
{\bf R}^\kappa= & \ R_{x^\kappa}+
\left[C_{\kappa}^{1}\,R_{x^{\kappa-1}u}-Q_{x^\kappa}\right]U^1+
\left[C_{\kappa}^{2}\,R_{x^{\kappa-2}u}-C_{\kappa}^{1}
Q_{x^{\kappa-1}}\right] U^2+\\
\ \ \ \ & \ +
\left[C_{\kappa}^{2}\,R_{x^2u} -C_{\kappa}^{3}\,Q_{x^3} 
\right] U^{\kappa-2}+
\left[C_{\kappa}^{1}\,R_{xu}-C_{\kappa}^{2}\,Q_{x^2}\right]
U^{\kappa-1}+\\
\ \ \ \ & \ +
\left[C_{\kappa}^{1}\,R_{u^2}-\kappa^2\,Q_{xu}\right]
U^1\,U^{\kappa-1}+
\left[-C_{\kappa+1}^{2}\,Q_u\right] U^2 U^{\kappa-1}+\\
\ \ \ \ & \ +
\left[R_u-C_{\kappa}^{1}\,Q_x\right] U^\kappa +
\left[-C_{\kappa+1}^{1}\,Q_u\right] U^1U^\kappa+\\
\ \ \ \ & \ + {\sf Remainder},
\endaligned\right.
\end{equation}
where the term {\sl Remainder} denotes the remaining terms in the
expansion of $R^\kappa$.
\end{lemma}
We note that the formula~(\ref{e210}) is valid for $\kappa=3$, comparing
with~(\ref{e29}), with the convention that the terms $U^{\kappa-2}$
and $U^{\kappa-1}$ vanish {\rm (}they coincide with $U^1$ and
$U^2${\rm )}, and replacing the coefficient $-C_{\kappa+1}^2\,
Q_u=-C_4^2\, Q_u=-6\, Q_u$ of the monomial $U^2\, U^{\kappa-1}$ by
$-3\, Q_u$, as it appears in~(\ref{e29}). The proof goes by a
straightforward computation, applying the recursive definition of this
partial formula.

\subsection{Computation of $R^\kappa$ in the general case} 
Following the exact same scheme as in the case $n=1$ we give 
the general partial formula for $R^\kappa$. We start with the
 first three families of coefficients ${\bf R}_{ k_1}^j$, ${\bf
R}_{k_1, k_2}^j$ and ${\bf R}_{ k_1, k_2, k_3}^j$. Let $\delta_p^q$ be
the Kronecker symbol, equal to $1$ if $p=q$ and to $0$ if $p\neq
q$. More generally, the {\sl generalized Kronecker symbols} are
defined by $ \delta_{ p_1, \dots, p_k}^{q_1, \dots, q_k}:= \delta_{
p_1}^{ q_1} \delta_{ p_2}^{ q_2}\cdots \delta_{ p_k}^{ q_k}$.

By convention, the indices $j$, $i_1$, $i_2$, $\dots$, $i_\lambda$ run
in the set $\{1,\dots,m\}$, the indices $k$, $k_1$, $k_2$, $\dots$,
$k_\lambda$ and $l$, $l_1$, $l_2$, $\dots$, $l_\lambda$ running in
$\{1,\dots,n\}$. Hence we will write $\sum_{i_1=1}^m\,
\sum_{i_2=1}^m\, \cdots \sum_{i_\lambda=1}^m$ as
$\sum_{i_1,\dots,i_\lambda}$ and $\sum_{l_1=1}^n\, \sum_{l_2=1}^n\,
\cdots \sum_{l_\lambda=1}^n$ as $\sum_{l_1,\dots,l_\lambda}$. The
letters $i_1,i_2,\dots,i_\lambda$ and $l_1,l_2,\dots,l_\lambda$ will
always be used for the summations in the development of ${\bf
R}_{k_1,k_2,\dots,k_{\lambda}}^j$.  We will always use the indices $j$
and $k_1,k_2,\dots, k_{\lambda}$ to write the coefficient ${\bf
R}_{k_1,k_2,\dots,k_{\lambda}}^j$.

We have:
\begin{equation}\label{e211}
\left\{
\aligned
{\bf R}_{k_1}^j=
& \
R_{x_{k_1}}^j+\sum_{i_1}\, 
\sum_{l_1}\, 
\left[
\delta_{k_1}^{l_1}\, R_{u^{i_1}}^j-
\delta_{i_1}^j \, Q_{x_{k_1}}^{l_1}
\right] U_{l_1}^{i_1}+\\
&
\ \ \ \ \ \ \ \ \ \ \ \ \ \ \ \ \ 
+\sum_{i_1,i_2}\, 
\sum_{l_1,l_2}\, 
\left[
-\delta_{i_2}^j\, \delta_{k_1}^{l_1} \, 
Q_{u^{i_1}}^{l_2}
\right] U_{l_1}^{i_1} \, U_{l_2}^{i_2}.
\endaligned\right.
\end{equation}
For ${\bf R}_{k_1,k_2}^j$ we have:
\begin{equation}\label{e212}
\left\{
\aligned
{}
&
{\bf R}_{k_1,k_2}^j
=
R_{x_{k_1}x_{k_2}}^j+
\sum_{i_1}\, \sum_{l_1}\, 
\left[
\delta_{k_2}^{l_1} \, 
R_{x_{k_1}u^{i_1}}^j+
\delta_{k_1}^{l_1}\, 
R_{x_{k_2}u^{i_1}}^j-
\delta_{i_1}^j\, Q_{x_{k_1}x_{k_2}}^{l_1}
\right]U_{l_1}^{i_1}+\\
& \
+\sum_{i_1,i_2}\, \sum_{l_1,l_2}\, 
\left[
\delta_{k_1,k_2}^{l_1,l_2}\, R_{u^{i_1}u^{i_2}}^j-
\delta_{i_2}^j \left(\delta_{k_1}^{l_1}\, 
Q_{x_{k_2}u^{i_1}}^{l_2}+ 
\delta_{k_2}^{l_1}\, Q_{x_{k_1}u^{i_1}}^{l_2}\right)
\right] 
U_{l_1}^{i_1} \, U_{l_2}^{i_2}+\\
& \ 
+
\sum_{i_1,i_2,i_3}\, \sum_{l_1,l_2,l_3}\,
\left[
-\delta_{i_3}^j\, \delta_{k_1,k_2}^{l_1,l_2}\, 
Q_{u^{i_1}u^{i_2}}^{l_3}
\right] 
U_{l_1}^{i_1} \, 
U_{l_2}^{i_2} \, U_{l_3}^{i_3}+ \\
& \
+
\sum_{i_1}\, \sum_{l_1,l_2}\, 
\left[
\delta_{k_1,k_2}^{l_1,l_2}\, R_{u^{i_1}}^j-
\delta_{i_1}^j\, \delta_{k_2}^{l_1}\, 
Q_{x_{k_1}}^{l_2}-
\delta_{i_1}^j\, \delta_{k_1}^{l_1} \, Q_{x_{k_2}}^{l_2}
\right]
U_{l_1,l_2}^{i_1}+ \\
& \
+
\sum_{i_1,i_2}\, \sum_{l_1,l_2,l_3}\, 
\left[
-\delta_{i_2}^j\, \delta_{k_1,k_2}^{l_1,l_2}\, 
Q_{u^{i_1}}^{l_3}-
\delta_{i_2}^j \, \delta_{k_1,k_2}^{l_3,l_1}\, 
Q_{u^{i_1}}^{l_2}-\delta_{i_1}^j\, 
\delta_{k_1,k_2}^{l_2,l_3}\, Q_{u^{i_2}}^{l_1}
\right]
U_{l_1}^{i_1}\, U_{l_2,l_3}^{i_2}.
\endaligned\right.
\end{equation}
Since we also treat systems of order $\kappa\geq 3$, 
it is necessary to compute ${\bf
R}_{k_1,k_2,k_3}^j$. We write this as follows:
\begin{equation}\label{e213}
{\bf R}_{k_1,k_2,k_3}^j={\rm I}+{\rm II}+{\rm III}, 
\end{equation}
where the first term {\rm I} involves only polynomials in $U_{l_1}^{i_1}$:
\begin{equation}\label{e214}
\left\{
\aligned
{\rm I}=
& \
R_{x_{k_1}x_{k_2}x_{k_3}}^j+
\sum_{i_1}\, \sum_{l_1}\, 
\left[
\delta_{k_1}^{l_1}\,R_{x_{k_2}x_{k_3}u^{i_1}}^j+
\delta_{k_2}^{l_1}\, R_{x_{k_1}x_{k_3}u^{i_1}}^j+
\delta_{k_3}^{l_1}\, R_{x_{k_1}x_{k_2}u^{i_1}}^j-\right. \\
& \
\left. -
\delta_{i_1}^j\, Q_{x_{k_1}x_{k_2}x_{k_3}}^{l_1}
\right]
U_{l_1}^{i_1}+
\sum_{i_1,i_2}\, \sum_{l_1,l_2}\, 
\left[
\delta_{k_1,k_2}^{l_1,l_2}\, R_{x_{k_3}u^{i_1}u^{i_2}}^j+
\delta_{k_3,k_1}^{l_1,l_2}\, R_{x_{k_2}u^{i_1}u^{i_2}}^j+ \right. \\
& \
\left. +
\delta_{k_2,k_3}^{l_1,l_2}\, R_{x_{k_1}u^{i_1}u^{i_2}}^j-
\delta_{i_2}^j\, \delta_{k_1}^{l_1}\, Q_{x_{k_2}x_{k_3}u^{i_1}}^{l_2}-
\delta_{i_2}^j\, \delta_{k_2}^{l_1}\, Q_{x_{k_1}x_{k_3}u^{i_1}}^{l_2}- 
\right. \\
& \
\left. -
\delta_{i_2}^j\, \delta_{k_3}^{l_1}\, Q_{x_{k_1}x_{k_2}u^{i_1}}^{l_2}
\right]
U_{l_1}^{i_1}\, U_{l_2}^{i_2}+
\sum_{i_1,i_2,i_3}\, \sum_{l_1,l_2,l_3}\, 
\left[
\delta_{k_1,k_2,k_3}^{l_1,l_2,l_3}\, R_{u^{i_1}u^{i_2}u^{i_3}}^j- \right. \\
& \
\left. -
\delta_{i_3}^j\, \delta_{k_1,k_2}^{l_1,l_2}\, Q_{x_{k_3}u^{i_1}u^{i_2}}^{l_3}-
\delta_{i_3}^j\, \delta_{k_2,k_3}^{l_1,l_2}\, Q_{x_{k_1}u^{i_1}u^{i_2}}^{l_3}-
\right. \\
& \ 
\left. -
\delta_{i_3}^j\, \delta_{k_1,k_3}^{l_1,l_2}\, Q_{x_{k_2}u^{i_1}u^{i_2}}^{l_3}
\right] 
U_{l_1}^{i_1} \, U_{l_2}^{i_2}\, 
U_{l_3}^{i_3}+ \\
& \
+\sum_{i_1,i_2,i_3,i_4}\, \sum_{l_1,l_2,l_3,l_4}\, 
\left[
-\delta_{i_4}^j\, \delta_{k_1,k_2,k_3}^{l_1,l_2,l_3}\, 
Q_{u^{i_1}u^{i_2}u^{i_3}}^{l_4}
\right] U_{l_1}^{i_1} \, U_{l_2}^{i_2} \, U_{l_3}^{i_3} \, 
U_{l_4}^{i_4},
\endaligned\right.
\end{equation}
the second term {\rm II} involves at least once the 
monomial $U_{l_1,l_2}^{i_1}$:

\begin{equation}\label{e215}
\left\{
\aligned
{\rm II}=
& \
\sum_{i_1}\, \sum_{l_1,l_2}\, 
\left[
\delta_{k_1,k_2}^{l_1,l_2}\, R_{x_{k_3}u^{i_1}}^j+
\delta_{k_3,k_1}^{l_1,l_2} \, R_{x_{k_2}u^{i_1}}^j+
\delta_{k_2,k_3}^{l_1,l_2}\, R_{x_{k_1}u^{i_1}}^j-\right. \\
& \
\left.
-\delta_{i_1}^j\left(
\delta_{k_1}^{l_1} \, Q_{x_{k_2}x_{k_3}}^{l_2}+
\delta_{k_2}^{l_1} \, Q_{x_{k_1}x_{k_3}}^{l_2}+
\delta_{k_3}^{l_1}\,  Q_{x_{k_1}x_{k_2}}^{l_2}
\right)
\right] U_{l_1,l_2}^{i_1}+\\
& \ +
\sum_{i_1,i_2}\, \sum_{l_1,l_2,l_3}\, 
\left[
\delta_{k_1,k_2,k_3}^{l_1,l_2,l_3}\, R_{u^{i_1}u^{i_2}}^j+
\delta_{k_1,k_2,k_3}^{l_3,l_1,l_2}\, R_{u^{i_1}u^{i_2}}^j+
\delta_{k_1,k_2,k_3}^{l_2,l_3,l_1}\, R_{u^{i_1}u^{i_2}}^j-\right. \\
& \
\left. -
\delta_{i_1}^j\left(
\delta_{k_1,k_2}^{l_2,l_3}\, Q_{x_{k_3}u^{i_2}}^{l_1}+
\delta_{k_3,k_1}^{l_2,l_3}\, Q_{x_{k_2}u^{i_2}}^{l_1}+
\delta_{k_2,k_3}^{l_2,l_3}\, Q_{x_{k_1}u^{i_2}}^{l_1}
\right)-\right.\\
& \ 
\left. -
\delta_{i_2}^j\left(
\delta_{k_1,k_2}^{l_1,l_2}\, Q_{x_{k_3}u^{i_1}}^{l_3}+
\delta_{k_3,k_1}^{l_1,l_2}\, Q_{x_{k_2}u^{i_1}}^{l_3}+
\delta_{k_2,k_3}^{l_1,l_2}\, Q_{x_{k_1}u^{i_1}}^{l_3}+\right.\right. \\
& \
\left. \left. +
\delta_{k_1,k_2}^{l_3,l_1}\, Q_{x_{k_3}u^{i_1}}^{l_2}+
\delta_{k_3,k_1}^{l_3,l_1}\, Q_{x_{k_2}u^{i_1}}^{l_2}+
\delta_{k_2,k_3}^{l_3,l_1}\, Q_{x_{k_1}u^{i_1}}^{l_2}
\right)
\right]U_{l_1}^{i_1}\, U_{l_2,l_3}^{i_2}+\\
& \ +
\sum_{i_1,i_2,i_3}\, \sum_{l_1,l_2,l_3,l_4}\, 
\left[
-\delta_{i_3}^j\left(
\delta_{k_1,k_2,k_3}^{l_1,l_2,l_3}\, Q_{u^{i_1}u^{i_2}}^{l_4}+
\delta_{k_1,k_2,k_3}^{l_1,l_4,l_2}\, Q_{u^{i_1}u^{i_2}}^{l_3}+\right.
\right. \\
& \ 
\left. \left.
\delta_{k_1,k_2,k_3}^{l_3,l_1,l_2}\, Q_{u^{i_1}u^{i_2}}^{l_4}
\right)- 
\delta_{i_1}^j\left(
\delta_{k_1,k_2,k_3}^{l_3,l_2,l_4}\, Q_{u^{i_2}u^{i_3}}^{l_1}+
\delta_{k_1,k_2,k_3}^{l_4,l_3,l_2}\, Q_{u^{i_2}u^{i_3}}^{l_1}+ \right. 
\right. \\
& \
\left. \left.
+\delta_{k_1,k_2,k_3}^{l_2,l_3,l_4}\, Q_{u^{i_1}u^{i_2}}^{l_1}
\right)
\right] U_{l_1}^{i_1} \, U_{l_2}^{i_2} \, U_{l_3,l_4}^{i_3}+
\sum_{i_1,i_2}\, \sum_{l_1,l_2,l_3,l_4}\, 
\left[
-\delta_{i_2}^j\left(
\delta_{k_1,k_2,k_3}^{l_1,l_2,l_3}\, Q_{u^{i_1}}^{l_4}+ \right. \right. \\
& \ 
\left. \left. +
\delta_{k_1,k_2,k_3}^{l_3,l_1,l_2}\, Q_{u^{i_1}}^{l_4}+
\delta_{k_1,k_2,k_3}^{l_2,l_3,l_1}\, Q_{u^{i_1}}^{l_4}
\right)
\right] U_{l_1,l_2}^{i_1}\, U_{l_3,l_4}^{i_2}
\endaligned\right.
\end{equation}
and the third term {\rm III} involves at least once the 
monomial $U_{l_1,l_2,l_3}^{i_1}$ (note that there is no 
term involving simultaneously $U_{l_1,l_2}^{i_1}$ and 
$U_{l_1,l_2,l_3}^{i_1}$):
\begin{equation}\label{e216}
\left\{
\aligned
{\rm III}=
& \
\sum_{i_1}\, \sum_{l_1,l_2,l_3}\, 
\left[
\delta_{k_1,k_2,k_3}^{l_1,l_2,l_3}\, 
R_{u^{i_1}}^j-\delta_{i_1}^j\left(
\delta_{k_2,k_3}^{l_1,l_2}\, Q_{x_{k_1}}^{l_3}+
\delta_{k_3,k_1}^{l_1,l_2}\, Q_{x_{k_2}}^{l_3}+ \right. \right. \\
& \
\left. \left. +
\delta_{k_1,k_2}^{l_1,l_2}\, Q_{x_{k_3}}^{l_3}
\right)
\right]U_{l_1,l_2,l_3}^{i_1}+
\sum_{i_1,i_2}\,\sum_{l_1,l_2,l_3,l_4}
\left[
-\delta_{i_1}^j\, \delta_{k_1,k_2,k_3}^{l_2,l_3,l_4}\, 
Q_{u^{i_2}}^{l_1}- \right. \\
& \ 
\left. -
\delta_{i_2}^j\left(
\delta_{k_1,k_2,k_3}^{l_1,l_2,l_3}\, Q_{u^{i_1}}^{l_4}+
\delta_{k_1,k_2,k_3}^{l_4,l_1,l_2}\, Q_{u^{i_1}}^{l_3}+
\delta_{k_1,k_2,k_3}^{l_3,l_4,l_1}\, Q_{u^{i_1}}^{l_2}
\right)
\right]U_{l_1}^{i_1}\, U_{l_2,l_3,l_4}^{i_2}.
\endaligned\right.
\end{equation}
Before giving the partial expression of $R^\kappa$ we introduce some
notations. For $p\in \N$ with $p\geq 1$, let $\mathfrak{S}_p$ be the
group of permutations of $\{1,2,\dots,p\}$.  For $q\in \N$ with $1\leq
q \leq p-1$, let $\mathfrak{S}_p^q$ be the set of permutations
$\sigma\in\mathfrak{S}_p$ such that $\sigma(1)< \sigma(2) < \cdots <
\sigma(q)$ and $\sigma(q+1) < \sigma(q+2) < \cdots < \sigma(p)$. Its
cardinal is $C_p^q$. Let $\mathfrak{C}_p$ be the group of cyclic
permutations of $\{1,2,\dots,p\}$. Reasoning recursively from the
formula of ${\bf R}_{k_1,k_2,k_3}^j$ given by~(\ref{e213}), we may
generalize Lemma~\ref{lem21}:

\begin{lemma}\label{lem22}
For every $\kappa\geq 4$ and for every $j=1,\dots,m$,
$k_1,\dots,k_\kappa=1,\dots,n$, we have:
\begin{equation}\label{e217}
{\bf R}_{k_1,k_2,\dots,k_\kappa}^j= I_1 + \cdots + I_9 + {\sf Remainder}
\end{equation}
where
\noindent $I_1 = R_{x_{k_1}x_{k_2}\dots x_{k_\kappa}}^j,$

\noindent $\aligned
I_2 = \sum_{i_1}\, \sum_{l_1}\, 
\left[
\sum_{\sigma\in\mathfrak{S}_\kappa^1}\, 
\delta_{k_{\sigma(1)}}^{l_1}\, 
R_{x_{k_{\sigma(2)}}\cdots x_{k_{\sigma(\kappa)}}u^{i_1}}^j-
\delta_{i_1}^j\, Q_{x_{k_1}\dots x_{k_\kappa}}^{l_1}
\right]U_{l_1}^{i_1},
\endaligned$

\noindent $\aligned
I_3=&\sum_{i_1}\, \sum_{l_1,l_2}\, 
\left[
\sum_{\sigma\in\mathfrak{S}_\kappa^2}\, 
\delta_{k_{\sigma(1)},k_{\sigma(2)}}^{l_1,l_2}\, 
R_{x_{k_{\sigma(3)}}\cdots x_{k_{\sigma(\kappa)}}u^{i_1}}^j-\right. \\
& \ 
\left. \ \ \ \ \ \ \ \ \ \ \ \ \ -  
\delta_{i_1}^j
\left(
\sum_{\sigma\in\mathfrak{S}_\kappa^1}\, 
\delta_{k_{\sigma(1)}}^{l_1}\, 
Q_{x_{k_{\sigma(2)}}\cdots x_{k_{\sigma(\kappa)}}}^{l_2}
\right)
\right] U_{l_1,l_2}^{i_1},
\endaligned$

\noindent $\aligned
I_4=&\sum_{i_1}\, \sum_{l_1,\dots,l_{\kappa-2}}\, 
\left[
\sum_{\sigma\in\mathfrak{S}_\kappa^{\kappa-2}}\, 
\delta_{k_{\sigma(1)},\dots,k_{\sigma(\kappa-2)}}^{l_1,\dots\dots,
l_{\kappa-2}}\,
R_{x_{k_{\sigma(\kappa-1)}}x_{k_{\sigma(\kappa)}}u^{i_1}}^j-\right. \\
& \
\left. \ \ \ \ \ \ \ \ \ \ \ \ \ -
\delta_{i_1}^j\left(
\sum_{\sigma\in\mathfrak{S}_\kappa^{\kappa-3}}\, 
\delta_{k_{\sigma(1)},\dots,k_{\sigma(\kappa-3)}}^{l_1,\dots\dots,
l_{\kappa-3}}\,
Q_{x_{k_{\sigma(\kappa-2)}}x_{k_{\sigma(\kappa-1)}}x_{k_{\sigma(\kappa)}}}^{
l_{\kappa-2}}
\right)
\right]U_{l_1,\dots,l_{\kappa-2}}^{i_1},
\endaligned$

\noindent $\aligned
I_5=&
\sum_{i_1}\, \sum_{l_1,\dots,l_{\kappa-1}}\, 
\left[
\sum_{\sigma\in\mathfrak{S}_\kappa^{\kappa-1}}\, 
\delta_{k_{\sigma(1)},\dots,k_{\sigma(\kappa-1)}}^{l_1,\dots\dots,
l_{\kappa-1}}\,
R_{x_{k_{\sigma(\kappa)}}u^{i_1}}^j- \right. \\
& \
\left. \ \ \ \ \ \ \ \ \ \ \ \ \ \ \ \ \ \ \ \ \ \ \ \ \ \ -
\delta_{i_1}^j\left(
\sum_{\sigma\in\mathfrak{S}_\kappa^{\kappa-2}}\, 
\delta_{k_{\sigma(1)},\dots,k_{\sigma(\kappa-2)}}^{l_1,\dots\dots,
l_{\kappa-2}}\,
_{x_{k_{\sigma(\kappa-1)}}x_{k_{\sigma(\kappa)}}}^{l_{\kappa-1}}
\right)
\right]U_{l_1,\dots,l_{\kappa-1}}^{i_1},
\endaligned$

\noindent $\aligned
I_6=&
\sum_{i_1,i_2}\, \sum_{l_1,\dots,l_\kappa}\,
\left[
\sum_{\tau\in\mathfrak{C}_\kappa}\, 
\delta_{k_1,\dots\dots,k_\kappa}^{l_{\tau(1)},\dots,l_{\tau(\kappa)}}\,
R_{u^{i_1}u^{i_2}}^j-
\delta_{i_1}^j\left(
\sum_{\sigma\in\mathfrak{S}_\kappa^{\kappa-1}}\, 
\delta_{k_{\sigma(1)},\dots,k_{\sigma(\kappa-1)}}^{l_1,\dots\dots,l_\kappa}\, 
Q_{x_{k_{\sigma(\kappa)}}u^{i_2}}^{l_1}\right) - \right. \\
& \ 
\left. -
\delta_{i_2}^j\left(
\sum_{\sigma\in\mathfrak{S}_\kappa^{\kappa-1}}\, 
\left(
\delta_{k_{\sigma(1)},\dots,k_{\sigma(\kappa-1)}}^{l_1,\dots\dots,
l_{\kappa-1}}\, 
Q_{x_{k_{\sigma(\kappa)}}u^{i_1}}^{l_\kappa}+\cdots + 
\delta_{k_{\sigma(1)},\dots,k_{\sigma(\kappa-1)}}^{l_3,\dots\dots,l_1}\, 
Q_{x_{k_{\sigma(\kappa)}}u^{i_2}}^{l_2}
\right)
\right)
\right]\times \\
 & \ \ \ \ \ \ \ \ \ \ \ \ \ \ \ \ \ \ \ \ \ \ \ \ \ \ \ \ \ \ \ 
\ \ \ \ \ \ \ \ \ \ \ \ \times  U_{l_1}^{i_1}\, U_{l_2,\dots,l_\kappa}^{i_2},
\endaligned$

\noindent $\aligned
I_7=&
\sum_{i_1,i_2}\,\sum_{l_3,\dots,l_{\kappa+1}}\, 
\left[
-\delta_{i_1}^j\left(
\delta_{k_1,\dots,k_\kappa}^{l_2,\dots,l_{\kappa+1}}\, 
Q_{u^{i_2}}^{l_1}+\cdots+
\delta_{k_1,\dots,k_\kappa}^{l_{\kappa+1},\dots,l_2}\, Q_{u^{i_2}}^{l_1}
\right) - \right. \\
& \ 
\left. \ \ \ \ \ \ \ \ \ \ \ \ \ \ \ \ \ \ \ \ \ \ \ \ \ -
\delta_{i_2}^j\left(
\sum_{\tau\in\mathfrak{S}_\kappa^2}\, 
\delta_{k_1,\dots\dots,k_\kappa}^{l_{\tau(1)},\dots,l_{\tau(\kappa)}}\, 
Q_{u^{i_1}}^{l_{\kappa+1}}
\right)
\right] U_{l_1,l_2}^{i_1}\, U_{l_3,\dots,l_{\kappa+1}}^{i_2},
\endaligned$

\noindent $\aligned
I_8=&
\sum_{i_1}\, \sum_{l_1,\dots,l_\kappa}\, 
\left[
\delta_{k_1,\dots,k_\kappa}^{l_1,\dots,l_\kappa}\, R_{u^{i_1}}^j-
\delta_{i_1}^j\left(
\sum_{\sigma\in\mathfrak{S}_\kappa^{\kappa-1}}\, 
\delta_{k_{\sigma(1)},\dots,k_{\sigma(\kappa-1)}}^{l_1,\dots\dots,
l_{\kappa-1}}\, 
Q_{x_{k_{\sigma(\kappa)}}}^{l_\kappa}
\right)
\right] U_{l_1,\dots,l_\kappa}^{i_1},
\endaligned$

\noindent $\aligned
I_9=&
\sum_{i_1,i_2}\, \sum_{l_1,\dots,l_{\kappa+1}}\, 
\left[
-\delta_{i_1}^j\, \delta_{k_1,\dots,k_\kappa}^{l_2,\dots,l_{\kappa+1}}\, 
Q_{u^{i_2}}^{l_1}-
\delta_{i_2}^j\left(
\delta_{k_1,\dots,k_\kappa}^{l_1,\dots,l_\kappa}\, Q_{u^{i_1}}^{l_{\kappa+1}}+
\cdots+
\delta_{k_1,\dots,k_\kappa}^{l_3,\dots,l_1}\, Q_{u^{i_1}}^{l_2}
\right)
\right]\times\\
 &\ \ \ \ \ \ \ \ \ \ \ \ \ \ \ \ \ \ \ \ \ \ \ \ \ \ \ \ \ \ \ \ 
\ \ \ \ \ \ \ \ \ \ \ \ \ \ \ \ \ \ \times U_{l_1}^{i_1}\, U_{l_2,\dots,
l_{\kappa+1}}^{i_2}
\endaligned$

\noindent
and where the term {\sf Remainder} denotes the remaining terms in 
the expansion of $R_{k_1,k_2,\dots,k_\kappa}^j$.
\end{lemma}
 In $I_6$ the summation on the upper indices $(l_1,\dots,l_\kappa)$
gets on all the circular permutations of $\{1,2,\dots,\kappa\}$ except
the identity. In $I_7$ the summation gets on all the circular
permutations of $\{2,3,\dots,\kappa+1\}$. In $I_9$ the summation gets
on all the circular permutations of $\{1,2,\dots,\kappa+1\}$ except
the one transforming $(l_1,l_2,\dots, l_{\kappa+1})$ into
$(l_2,l_3,\dots,l_1)$.  For $\kappa=3$, comparing with~(\ref{e213}),
we see that the formula remains valid, with the same conventions as in
the case $n=1$.

\subsection{Lie criterion and defining equations of
$\mathfrak{Sym}(\mathcal{E})$} We recall the Lie criterion, presented
in Subsection~2.6 ({\it see} Theorem~2.71 of~\cite{ol1986}):

{\it A vector field $X$ is an infinitesimal symmetry of the completely
integrable system $(\mathcal E)$ if and only if its prolongation
$X^{(\kappa)}$ of order $\kappa$ is tangent to the skeleton
$\Delta_\mathcal E$ in the jet space $\mathcal{J}_{n,m}^\kappa$}.

The set of infinitesimal symmetries of $(\mathcal{E})$ forms a Lie
algebra, since we have the relation $[X,X']^{(\kappa)}=[X^{(\kappa)},
{X'}^{(\kappa)}]$ ({\it cf.}~\cite{ol1986}). We will denote by
$\mathfrak{Sym}(\mathcal{E})$ this Lie algebra.  The aim of the
forecoming Section is to obtain precise bounds on the dimension of the
Lie algebra $\mathfrak{Sym}(\mathcal{E})$ of infinitesimal symmetries
of ($\mathcal E$). For simplicity we start with the case $n=m=1$.

\section{Optimal upper bound on $\dim_\K\mathfrak{Sym}(\mathcal{E})$ 
when $n=m=1$.} 

\subsection{Defining equations for $\mathfrak{Sym}(\mathcal{E})$.} 
Applying the Lie criterion, the tangency condition of $X^{(\kappa)}$ 
to $\Delta_\mathcal{E}$ is equivalent to the identity:
\begin{equation}\label{e31}
{\bf R}^\kappa-\left[ 
Q \, \frac{\partial F}{\partial x}+R \, 
\frac{\partial F}{\partial u}+
{\bf R}^1 \, \frac{\partial F}{\partial U^1} + {\bf R}^2 
\, \frac{\partial F}{\partial
U^2}+\cdots+{\bf R}^{\kappa-1} \, 
\frac{\partial F}{\partial U^{\kappa-1}}
\right]\equiv 0,
\end{equation} 
on the subvariety $\Delta_\mathcal{ E}$, that is to a formal identity
in $\K\{x, u, U^1, \dots, U^{\kappa-1}\}$, in which we replace the
variable $U^\kappa$ by $F(x,u, U^1, \dots, U^{ \kappa-1})$ in the two
monomials $U^\kappa$ and $U^1\, U^\kappa$ of ${\bf R}^\kappa$, {\it
cf.}  Lemma~\ref{lem21}. Expanding $F$ and its partial derivatives in
power series of the variables $(U^1,\dots,U^{\kappa-1})$ with analytic
coefficients in $(x,u)$, we may rewrite~(\ref{e31}) as follows:
\begin{equation}\label{e32}
\left\{
\aligned
\sum_{\mu_1,\dots,\mu_{\kappa-1}\geq 0}\,
\left[
\Phi_{\mu_1,\dots,\mu_{\kappa-1}}
\left(x,u,(Q_{x^ku^l})_{k+l\leq \kappa}, 
(R_{x^ku^l})_{k+l\leq \kappa}\right)\right]\times \\ 
\times (U^1)^{\mu_1}\dots (U^{\kappa-1})^{\mu_{\kappa-1}}\equiv 0,
\endaligned\right.
\end{equation}
where the expressions 
\begin{equation}\label{e33}
\Phi_{\mu_1,\dots,\mu_{\kappa-1}}
\left(x,u,(Q_{x^ku^l})_{k+l\leq \kappa}, 
(R_{x^ku^l})_{k+l\leq \kappa}\right)
\end{equation}
are {\sl linear}\, with respect to the partial derivatives
$((Q_{x^ku^l})_{k+l\leq \kappa}, (R_{x^ku^l})_{k+l\leq \kappa})$, with
analytic coefficients in $(x,u)$.  By construction these coefficients
essentially depend on the expansion of $F$. The tangency
condition~(\ref{e32}) is equivalent to the following infinite linear
system of partial differential equations, called {\sl defining
equations of $\mathfrak{ Sym} (\mathcal{E})$}:
\begin{equation}\label{e34}
\Phi_{\mu_1,\dots,\mu_{\kappa-1}}
\left(x,u,(Q_{x^ku^l}(x,u))_{k+l\leq \kappa}, 
(R_{x^ku^l}(x,u))_{k+l\leq \kappa}\right)=0, 
\end{equation}
satisfied by $(Q(x,u), R(x,u))$. The Lie method consists in studying
the solutions of this linear system of partial differential equations.

\subsection{Homogeneous system}\label{sub3.1}

As mentioned in the introduction, we focus our attention on the case
$\kappa \geq 3$. Denote by $(\mathcal{E}_0)$ the homogeneous equation
$u_{x^\kappa}=0$ of order $\kappa$. The general solution
$u=\sum_{l=0}^{\kappa-1}\, \lambda_l\, x^l$ consists of polynomials of
degree $\leq \kappa-1$ and the defining equation~(\ref{e31}) reduces
to ${\bf R}^\kappa=0$. Using the expression~(\ref{e210}),
expanding~(\ref{e32}), (\ref{e33}) and considering only the
coefficients of the five monomials ${\rm ct.}$, $U^{\kappa-2}$,
$U^{\kappa-1}$, $U^1\, U^{\kappa-1}$ and $U^2\, U^{\kappa-1}$, we
obtain the five following partial differential equations, which are
sufficient to determine $\mathfrak{Sym}(\mathcal{E}_0)$:
\begin{equation}\label{e35}
\left\{
\aligned
R_{x^\kappa} =
& \
0, \\
R_{x^2u}-\frac{(\kappa-2)}{3} \, 
Q_{x^3}=
& \
0,\\
R_{xu}-\frac{(\kappa-1)}{2} \, Q_{x^2}=
& \ 
0,\\
R_{u^2}-\kappa\, Q_{xu}= 
& \
0,\\
Q_u=
& \
0.
\endaligned\right.
\end{equation}
The general solution of this system is evidently:
\begin{equation}\label{e36}
\left\{
\aligned
Q= & \ A+B\,x+C\,x^2,\\
R= & \ (\kappa-1)\, 
C\, xu+D\,u+E^0+E^1\, x+\cdots+
E^{\kappa-1}\, x^{\kappa-1},
\endaligned\right.
\end{equation}
where the $(\kappa+4)$ constants $A,\, B,\, C,\, D,\, E^0, \, E^1,
\dots, \, E^{\kappa-1}$ are arbitrary. Computing explicitely the flows
of the $(\kappa+4)$ generators $\partial/\partial x$, $x\partial
/\partial x$, $x^2\, \partial /\partial x+(\kappa-1)\, xu\,
\partial/\partial u$, $u\,\partial /\partial u$, $\partial /\partial
u$, $x\, \partial /\partial u,\dots, x^{\kappa-1}\, \partial/\partial
u$, we check easily that they stabilize the graphs of polynomials of
degree $\leq \kappa-1$. Moreover they span a Lie algebra of dimension
$(\kappa +4)$ and the general form of a Lie symmetry is:
\begin{equation}\label{e37}
(x,u)\longmapsto 
\left(
\frac{\alpha_0+\alpha_1 x}{1+\varepsilon x}, \ \
\frac{\beta u+\gamma_0+\gamma_1x+\cdots+\gamma_{\kappa-1} 
x^{\kappa-1}}{(1+\varepsilon x)^{\kappa-1}}
\right).
\end{equation}

\subsection{Nonhomogeneous system}
Consider for $\kappa \geq 3$ the equation~(\ref{e31}) after replacing
the variable $U^\kappa$ by $F$. Let $\Phi(U^\lambda)$ denote an
arbitrary term of the form $\phi(x,u)\, U^\lambda$, where $\phi(x,u)$
is an analytic function. We consider the five following terms
$\Phi({\rm ct.})$, $\Phi(U^{\kappa-2})$, $\Phi(U^{\kappa-1})$,
$\Phi(U^1\, U^{\kappa-1})$ and $\Phi(U^2\, U^{\kappa-1})$. Since some
multiplications of monomials appear in the expression~(\ref{e31}), we
must be aware of the fact that $\Phi(U^1\, U^{\kappa-1})\equiv
\Phi(U^1)\, \Phi(U^{\kappa-1})$ and $\Phi(U^2\, U^{\kappa-1})\equiv
\Phi(U^2)\, \Phi(U^{\kappa-1})$. Consequently in the expansion
of~(\ref{e31}) we must take into account the seven types of monomials
$\Phi({\rm ct.})$, $\Phi(U^1)$, $\Phi(U^2)$, $\Phi(U^{\kappa-2})$,
$\Phi(U^{\kappa-1})$, $\Phi(U^1\, U^{\kappa-1})$ and $\Phi(U^2\,
U^{\kappa-1})$.  The $(\kappa+1)$ derivatives $\partial F/\partial x$,
$\partial F/\partial u$, $\partial F/\partial U^1,\dots,\partial
F/\partial U^{\kappa-1}$ appearing in the brackets of~(\ref{e31}), and
the term $F$ appearing in the expression of ${\bf R}^\kappa$ after
replacing $U^\kappa$ by $F$ ({\it cf.} the last two monomials
$U^\kappa$ and $U^1\, U^\kappa$ in~(\ref{e210})) may all contain the
seven monomials ${\rm ct.}$, $U^1$, $U^2$, $U^{\kappa-2}$,
$U^{\kappa-1}$, $U^1\, U^{\kappa-1}$ and $U^2\, U^{\kappa-1}$. For $F$
and its $(\kappa+1)$ first derivatives we use the generic simplified
notation
\begin{equation}\label{e38}
\Phi({\rm ct.})+\Phi(U^1)+\Phi(U^2)+\Phi(U^{\kappa-2})+
\Phi(U^{\kappa-1})+\Phi(U^1\, 
U^{\kappa-1})+
\Phi(U^2\, U^{\kappa-1}),
\end{equation} 
to name the seven monomials appearing {\it a priori}.  Hence,
expanding~(\ref{e31}), picking up the only terms which may contain the
five monomials we are interested in, and using the formula of
Lemma~\ref{lem21} for ${\bf R}^\lambda$ $(1 \leq \lambda \leq
\kappa)$, we obtain the following expression:
\begin{equation}\label{e39}
\left\{
\aligned
{}
&
R_{x^\kappa}+\left[
C_{\kappa}^{2}\,
R_{x^2u}-C_{
\kappa}^{3}\, 
Q_{x^3}\right]
U^{\kappa-2}+
\left[
C_{\kappa}^{1} \, R_{xu}-
C_{\kappa}^{2}\, Q_{x^2}
\right] U^{\kappa-1}+\\
& \
+\left[C_{\kappa}^{1}\, 
R_{u^2}-\kappa^2\, Q_{xu} \right] 
U^1\, U^{\kappa-1}+
\left[-C_{\kappa+1}^{2}\, Q_u\right]U^2\, U^{\kappa-1}+\\
& \
+\left\{
R_u-C_{\kappa}^{1} \, Q_x +
\left[
-C_{\kappa+1}^{1}\, Q_u
\right] U^1
\right\}\times
\\
& \
\times\left\{
\Phi({\rm ct.})+\Phi(U^1)+\Phi(U^2)+\Phi(U^{\kappa-2})+
\Phi(U^{\kappa-1})+
\Phi(U^1\, U^{\kappa-1})+
\Phi(U^2\, U^{\kappa-1})
\right\}-\\
& \
-\left\{
Q+R+R_x+\left[R_u-Q_x\right]U^1+R_{x^2}+
\left[2R_{xu}-Q_{x^2}\right]U^1+\right.\\
& \
\left.
+\left[R_u-2Q_x\right]U^2+
\dots+R_{x^{\kappa-3}}+
\left[C_{\kappa-3}^1\, R_{x^{\kappa-4}u}-
Q_{x^{\kappa-3}}\right]U^1+\right.\\
& \
\left.
+\left[
C_{\kappa-3}^2\, R_{x^{\kappa-5}u}-C_{\kappa-3}^1\, 
Q_{x^{\kappa-4}}\right] U^2
+R_{x^{\kappa-2}}+\right.\\
& \
\left.
+
\left[
C_{\kappa-2}^{1}\, R_{x^{\kappa-3}u}-Q_{x^{\kappa-2}}
\right] U^1+
\left[
C_{\kappa-2}^2\, R_{x^{\kappa-4}u}-
C_{\kappa-2}^1\, Q_{x^{\kappa-3}}
\right] U^2+\right.\\
& \
\left.
+\left[
R_u-C_{\kappa-2}^1\, Q_x
\right] U^{\kappa-2}+
\left[-C_{\kappa-1}^1\, Q_u
\right] U^1\, U^{\kappa-2}+\right.\\
& \
\left.
+
R_{x^{\kappa-1}}+
\left[
C_{\kappa-1}^1\, R_{x^{\kappa-2}u}-
Q_{x^{\kappa-1}}
\right] U^1+\right.\\
& \
\left.
+\left[
C_{\kappa-1}^2\, R_{x^{\kappa-3}u}-
C_{\kappa-1}^1\, Q_{x^{\kappa-2}}
\right] U^2+
\left[C_{\kappa-1}^1\, R_{xu}-
C_{\kappa-1}^2\, 
Q_{x^2}\right] U^{\kappa-2}+\right.\\
& \
\left.
+\left[C_{\kappa-1}^1\, R_{u^2}-(\kappa-1)^2\, 
Q_{xu}\right] U^1\, U^{\kappa-2}
+\left[
R_u-C_{\kappa-1}^1\, Q_x
\right] U^{\kappa-1}+\right.\\
& \
\left.
+\left[-C_\kappa^1\, Q_u\right]
U^1\, U^{\kappa-1}
\right\}
\times\\
& \
\times\left\{
\Phi({\rm ct.})+\Phi(U^1)+\Phi(U^2)+\Phi(U^{\kappa-2})+
\Phi(U^{\kappa-1})+
\Phi(U^1\, U^{\kappa-1})+
\Phi(U^2\, U^{\kappa-1})
\right\}\\
& \
+ {\sf Remainder} \ 
\equiv 0.
\endaligned\right.
\end{equation}
Here the term {\sf Remainder} consists of the monomials, in the jet
variables, different from the five ones we are concerned with. The
first four lines before the sign ``$-$'' develop ${\bf R}^\kappa$ and
the third line consists of the factor $F$ replaced by~(\ref{e38}). In
the last line (note that this is multiplied by the nine preceding
lines) we replaced the $(\kappa+1)$ first partial derivatives of $F$
appearing in~(\ref{e31}) by the term~(\ref{e38}) which we factorized. 

By expanding the product appearing in this expression~(\ref{e39}), and
equaling to zero the coefficients of the five monomials ${\rm ct.}$,
$U^{\kappa-2}$, $U^{\kappa-1}$, $U^1\, U^{\kappa-1}$ and $U^2\,
U^{\kappa-1}$, we obtain the five following partial differential
equations
\begin{equation}\label{e310}
\left\{
\aligned
R_{x^\kappa}= & \ \Pi(x,u,Q,Q_x,R,R_x,\dots,R_{x^{\kappa-1}},R_u),\\
C_\kappa^2\,R_{x^2u}-C_\kappa^3 \, Q_{x^3}= &\ \Pi(x,u,Q,Q_x,Q_{x^2},R,R_x,
\dots,R_{x^{\kappa-1}},R_u,R_{xu}),\\
C_\kappa^1\,R_{xu}-C_\kappa^2 \, Q_{x^2}= & \ \Pi(x,u,Q,Q_x,R,R_x,\dots,
R_{x^{\kappa-1}},R_u),\\
C_\kappa^1\,R_{u^2}-\kappa^2 \, Q_{xu}= & \ \Pi(x,u,Q,Q_x,\dots,
Q_{x^{\kappa-1}},Q_u,R,R_x,\dots R_{x^{\kappa-1}},\\
& \ \ \ \ \ \ \ \ \ \ \ \ \ \ \ \ \ \ \ \ \ \ \ \ \ \ \ \ \ \ \ \ \ \ \ \ \ \
\ \ \ \ \ \ 
R_u,R_{xu},\dots,R_{x^{\kappa-2}u}),\\
-C_{\kappa+2}^2\,Q_u= & \ \Pi(x,u,Q,Q_x,\dots,
Q_{x^{\kappa-2}},R,R_x,\dots R_{x^{\kappa-1}},\\
& \ \ \ \ \ \ \ \ \ \ \ \ \ \ \ \ \ \ \ \ \ \ \ \ \ \ \ \ \ \ \ \ \ \ \ \ \ \
R_u,R_{xu},\dots,R_{x^{\kappa-3}u}).\\
\endaligned\right.
\end{equation}
Here by convention $\Pi$ denotes any linear quantity in $Q$, $R$ 
and some of their derivatives, of the form
\begin{equation}\label{eee}
\left\{
\aligned
{}
&
\Pi(x,u,Q_{x^{a_1}u^{b_1}}, \dots,
Q_{x^{a_p}u^{b_p}},R_{x^{c_1}u^{d_1}}, \dots,
R_{x^{c_q}u^{d_q}})=\\
& \ \ \ \ 
=\sum_{i=1}^p\, 
\phi_i(x,u)\, Q_{x^{a_i}u^{b_i}}(x,u)+
\sum_{j=1}^q \, 
\psi_j(x,u)\, R_{x^{c_j}u^{d_j}}(x,u),
\endaligned\right.
\end{equation}
where $\phi_i$ and $\psi_j$ are analytic in $(x,u)$. For instance, the
differentiation of $\Pi(x,u,Q,R,R_u)$ with respect to $x$ gives the
expression $\Pi(x,u,Q,Q_x,R,R_x,R_{xu})$. Let us introduce the
following collection of $(\kappa +4)$ partial derivatives of $(Q, R)$
defined by $J:=(Q, Q_x, Q_{x^2}, R, R_x, \dots, R_{x^{ \kappa-1 }},
R_u)$. The aim is now to make linear substitutions on the
system~(\ref{e310}) to obtain the system~(\ref{e319}) where the five
second members depend only on the collection $J$. The desired estimate
$\dim_\K \, \mathfrak{ Sym} (\mathcal{ E})\leq \kappa+4$ will follow
from~(\ref{e319}).

Let us differentiate the third equation
of~(\ref{e310}) with respect to $x$. Dividing by $C_\kappa^1$ we obtain:
\label{e311}
\begin{equation}
R_{x^2u}-\frac{(\kappa-1)}{2} \, Q_{x^3}=\Pi(x,u,Q,Q_x,Q_{x^2},R,R_x,\dots,
R_{x^\kappa},R_u,R_{xu}).
\end{equation}
Solving $R_{x^2u}$ and $Q_{x^3}$ by the second equality in~(\ref{e310}) 
and by~(\ref{e311}) we find
\begin{equation}\label{e312}
\left\{
\aligned
Q_{x^3}= & \ \Pi(x,u,Q,Q_x,Q_{x^2},R,R_x,\dots,R_{x^\kappa},R_u,R_{xu}),\\
R_{x^2u}= & \ \Pi(x,u,Q,Q_x,Q_{x^2},R,R_x,\dots,R_{x^\kappa},R_u,R_{xu}).
\endaligned\right.
\end{equation}
Replacing $R_{x^\kappa}$ by its value given by the first equality 
in~(\ref{e310}) we obtain for $Q_{x^3}$:
\begin{equation}\label{e313}
Q_{x^3}=\Pi(x,u,Q,Q_x,Q_{x^2},R,R_x,\dots,R_{x^{\kappa-1}},R_u,R_{xu}).
\end{equation}
If we write the third equality in~(\ref{e310}) as
\begin{equation}\label{e314}
R_{xu}=\Pi(x,u,Q,Q_x,Q_{x^2},R,R_x,\dots,R_{x^{\kappa-1}},R_u),
\end{equation}
we may replace $R_{xu}$ in~(\ref{e313}). This gives the desired 
dependence of $Q_{x^3}$ on the collection $J$:
\begin{equation}\label{e315}
Q_{x^3}=\Pi(x,u,Q,Q_x,Q_{x^2},R,R_x,\dots,R_{x^{\kappa-1}},R_u).
\end{equation}
We may now differentiate the equalities~(\ref{e314}) and~(\ref{e315})
with respect to $x$ up to the order $l$. At each differentiation we
replace $Q_{x^3}$, $R_{xu}$ and $R_{x^\kappa}$ by their values
in~(\ref{e314}), in (\ref{e315}) and in the first equality
in~(\ref{e310}) respectively. We obtain for $l\in\N$:
\begin{equation}\label{e316}
\left\{
\aligned
Q_{x^{l+3}}= & \ \Pi(x,u,Q,Q_x,Q_{x^2},R,R_x,\dots,R_{x^{\kappa-1}},R_u), \\
R_{x^{l+1}u}= & \ \Pi(x,u,Q,Q_x,Q_{x^2},R,R_x,\dots,R_{x^{\kappa-1}},R_u). 
\endaligned\right.
\end{equation}
Replacing these values in the fifth equality of~(\ref{e310}), we obtain
\begin{equation}\label{e317}
Q_u=\Pi(x,u,Q,Q_x,Q_{x^2},R,R_x,\dots,R_{x^{\kappa-1}},R_u).
\end{equation}
By replacing the fourth equality of~(\ref{e310}) we obtain finally 
\begin{equation}\label{e318}
R_{u^2}=\Pi(x,u,Q,Q_x,Q_{x^2},R,R_x,\dots,R_{x^{\kappa-1}},R_u).
\end{equation}
To summarize, using the first equality of~(\ref{e310}), using (\ref{e317}), 
(\ref{e318}), (\ref{e314}) and ~(\ref{e315}), we obtained the desired system:
\begin{equation}\label{e319}
\left\{
\aligned
R_{x^\kappa}=\Pi(x,u,Q,Q_x,Q_{x^2},R,R_x,\dots,R_{x^{\kappa-1}},R_u),\\
Q_u =\Pi(x,u,Q,Q_x,Q_{x^2},R,R_x,\dots,R_{x^{\kappa-1}},R_u),\\
R_{u^2}=\Pi(x,u,Q,Q_x,Q_{x^2},R,R_x,\dots,R_{x^{\kappa-1}},R_u),\\
R_{xu}=\Pi(x,u,Q,Q_x,Q_{x^2},R,R_x,\dots,R_{x^{\kappa-1}},R_u),\\
Q_{x^3}=\Pi(x,u,Q,Q_x,Q_{x^2},R,R_x,\dots,R_{x^{\kappa-1}},R_u).\\
\endaligned\right.
\end{equation}
We recall that the terms $\Pi$ are linear expressions of the
form~(\ref{eee}). Let us differentiate every equation of
system~(\ref{e319}) with respect to $x$ at an arbitrary order and let
us replace in the right hand side the terms $R_{x^{\kappa}}$, $R_{xu}$
and $Q_{x^3}$ that may appear at each step by their value
in~(\ref{e319}), and then differentiate with respect to $u$ at an
arbitrary order. We deduce that {\it all the partial derivatives of
the five functions $R_{ x^\kappa}$, $Q_u$, $R_{u^2}$, $R_{xu}$ and
$Q_{x^3}$ are also linear functions of the $(\kappa +4)$ partial
derivatives $(Q, Q_x, Q_{x^2}, R, R_x, \dots, R_{x^{ \kappa-1}},
R_u)$}. Thus the analytic functions $Q$ and $R$ are determined
uniquely by the value at the origin of the $(\kappa +4)$ partial
derivatives $(Q, Q_x, Q_{x^2}, R, R_x, \dots, R_{x^{ \kappa-1}},
R_u)$. This ends the proof of the inequality $\dim_\K \, \mathfrak{
Sym}\, (\mathcal{ E})\leq \kappa+4$. \qed

\section{Optimal upper bound on $\dim_\K \, \mathfrak{Sym}(\mathcal{E})$ 
in the general dimensional case}

\subsection{Defining equations for
$\mathfrak{Sym}(\mathcal{E})$} In the general dimensional case, the
tangency condition of the prolongation $X^\kappa$ of $X$ to the
skeleton gives the following equations for $j=1,\dots,m$ and
$k_1,\dots,k_\kappa=1,\dots,n$:
\begin{equation}\label{e41}
\left\{
\aligned
{}
&
{\bf R}_{k_1,\dots,k_\kappa}^j-\left[
\sum_{l=1}^n\, Q^l\, 
\frac{\partial F_{k_1,\dots,k_\kappa}^j}{\partial x_l}+
\sum_{i=1}^m\, R^i \,
\frac{\partial F_{k_1,\dots,k_\kappa}^j}{\partial u^i}+ \right. \\
& \ 
\left. +
\sum_{i_1}\, \sum_{l_1}\, 
{\bf R}_{l_1}^{i_1}\, 
\frac{\partial F_{k_1,\dots,k_\kappa}^j}{\partial U_{l_1}^{i_1}}+
\cdots+
\sum_{i_1}\, \sum_{l_1,\dots,l_{\kappa-1}}\, 
{\bf R}_{l_1,\dots,l_{\kappa-1}}^{i_1}\,
\frac{\partial F_{k_1,\dots,k_\kappa}^j}{
\partial U_{l_1,\dots,l_{\kappa-1}}^{i_1}}
\right]\equiv 0,
\endaligned\right.
\end{equation}
on $\Delta_\mathcal{E}$, by replacing the variables
$U_{l_1,\dots,l_\kappa}^{i_1}$ by $F_{l_1,\dots,l_\kappa}^{i_1}$
wherever they appear.  Let us expand $F_{k_1,\dots,k_\kappa}^j$ and
their partial derivatives and use the fact that ${\bf
R}_{k_1,\dots,k_\lambda}^j$ are polynomials expressions of the jets
variables $(U_{l_1}^{i_1},\dots,U_{l_1,\dots,l_\lambda}^{i_1})$, with
coefficients being linear expressions of the partial derivatives of
order $\leq \lambda+1$ of $Q^l$ and $R^j$. We obtain for $j=1,\dots,m$
and $k_1,\dots,k_\kappa=1,\dots,n$ some identities of the form
\begin{equation}\label{e42}
{\small 
\left\{
\aligned
\sum_{i_1, \ \dots, l_1,\dots}
\Phi_{k_1,\dots,k_\kappa; \ l_1,\dots\dots}^{j;i_1,\dots\dots}
\left(x,u,
(Q_{x^\alpha u^\beta}^l)_{1\leq l\leq n, \, \vert \alpha \vert+
\vert \beta \vert
\leq \kappa+1},
(R_{x^\alpha u^\beta}^j)_{1\leq j\leq m, \, \vert \alpha \vert+
\vert \beta \vert
\leq \kappa+1}
\right)\times
\\
\times U_{l_1}^{i_1}\dots U_{l_{\mu_1}}^{i_{\mu_1}}\times
U_{l_{\mu_1}+1,l_{\mu_1}+2}^{i_{\mu_1}+1}\cdots 
U_{l_{\mu_1+2\mu_2}-1}^{i_{\mu_1+\mu_2-1}}\, 
U_{l_{\mu_1+2\mu_2}}^{i_{\mu_1+\mu_2}}\times \cdots \cdots
\equiv 0,
\endaligned\right.
}
\end{equation}

\noindent satisfied if and only if the functions $Q^l$ and $R^j$ 
are solutions of the following system of partial differential equations
\begin{equation}\label{e43}
\Phi_{k_1,\dots,k_\kappa; \ l_1,\dots\dots}^{j,i_1,\dots\dots}
\left(x,u,
(Q_{x^\alpha u^\beta}^l)_{1\leq l\leq n, \, \vert \alpha \vert+
\vert \beta \vert
\leq \kappa+1},
(R_{x^\alpha u^\beta}^j)_{1\leq j\leq m, \, \vert \alpha \vert+
\vert \beta \vert
\leq \kappa+1}
\right)=0.
\end{equation}

\subsection{Homogeneous system}
We start by giving the general form of the symmetries of the
homogeneous system in the case $\kappa =2$. Then we prove the equality
$\dim_\K(\mathfrak{Sym}(\mathcal E_0)) = n^2+2n+m^2+ m \, C_{n+\kappa
-1}^{\kappa -1}$ in the case $\kappa \geq 3$.

\vskip 0,1cm
\noindent In the case $\kappa=2$ we obtain:
\begin{equation}\label{e44}
\left\{
\aligned
Q^l(x,u)=
& \
A^l+
\sum_{k_1=1}^n\, 
B_{k_1}^l\, x_{k_1}+
\sum_{i_1=1}^m\, C_{i_1}^l\, u^{i_1}+ \\
& \ \ \ \ \ \ \ \ \ \ \ \ \ \ \ \ \ \ +
\sum_{k_1=1}^n\, 
D_{k_1}\, x_l \, x_{k_1} +
\sum_{i_1=1}^m\, E_{i_1}\, x_l \, u^{i_1},\\
R^j(x,u)=
& \
F^j+\sum_{k_1=1}^n\, G_{k_1}^j\, x_{k_1}+
\sum_{i_1=1}^m \, H_{i_1}^j\, u^{i_1}+\\
& \ \ \ \ \ \ \ \ \ \ \ \ \ \ \ \ \ \ +
\sum_{k_1=1}^n\, D_{k_1}\, x_{k_1}\, u^j+
\sum_{i_1=1}^m\, E_{i_1}\, u^{i_1}\, u^j.
\endaligned\right.
\end{equation}
Here the $(n+m)(n+m+2)$ constants $A^ l,\, B_{k_1}^l, \, C_{i_1}^l, \,
D_{k_1}, \, E_{i_1}, \, F^j, \, G_{k_1}^j, \, H_{i_1}^j\in\K$ are
arbitrary.  Moreover one can check that the vector space spanned by
the $(n+m)(n+m+2)$ vector fields
\begin{equation}\label{e45}
\left\{
\aligned
{}
&
\frac{\partial }{\partial x_{k_1}}, \ x_{k_1}\, 
\frac{\partial }{\partial x_{k_2}}, \, 
u^{i_1}\, 
\frac{\partial }{\partial x_{k_1}},\\
&
x_{k_1}\left(x_1 \, \frac{\partial }{\partial x_1}+\cdots+
x_n\, \frac{\partial }{\partial x_n}+
u^1\, \frac{\partial }{\partial u^1}+\cdots+
u^m\, \frac{\partial }{\partial u^m}\right),\\
&
u^{i_1}\left(x_1 \, \frac{\partial }{\partial x_1}+\cdots+
x_n\, \frac{\partial }{\partial x_n}+
u^1\, \frac{\partial }{\partial u^1}+\cdots+
u^m\, \frac{\partial }{\partial u^m}\right),\\
&
\frac{\partial }{\partial u^{i_1}}, \ \
x_{k_1}\, \frac{\partial }{\partial u^{i_1}}, \ \
u^{i_1}\, \frac{\partial }{\partial u^{i_2}}
\endaligned\right.
\end{equation}
is stable under the Lie bracket action and that the flow of each of
these generators is a Lie symmetry of the system
$(\mathcal{E}_0)$. This proves that $\mathfrak{Sym}(\mathcal{E}_0)$ is
indeed a {\it Lie algebra}\, with dimension $(n+m)(n+m+2)$. Finally
the corresponding transformations close to the identity mapping are
projective, represented by the formula:
\begin{equation}\label{e46}
\left\{
\aligned
(x,u)\longmapsto 
\left(
\left(
\frac{\alpha_{l,0}+\sum_{k=1}^n\, 
\alpha_{l,k}\, x_k+\sum_{i=1}^m\, \alpha_{l,n+i}\, u^i}{
1+\sum_{k=1}^n\, 
\gamma_k\, x_k+\sum_{i=1}^m\, \gamma_{n+i}\, u^i}
\right)_{1\leq l\leq n}, \right. \\
\ \ \ \ \ \ \ \ \ \ \ 
\left.
\left(
\frac{\beta_{j,0}+\sum_{k=1}^n\, 
\beta_{j,k}\, x_k+\sum_{i=1}^m\, \beta_{j,n+i}\, u^i}{
1+\sum_{k=1}^n\, 
\gamma_k\, x_k+\sum_{i=1}^m\, \gamma_{n+i}\, u^i}
\right)_{1\leq j \leq m}
\right).
\endaligned\right.
\end{equation}
It is clear that these transformations preserve all the solutions of
$(\mathcal{E}_0): \ \ u_{x_{k_1}x_{k_2}}^j=0$, the graphs of affine
maps from $\K^n$ to $\K^m$.

In the case $\kappa \geq 3$ we consider the homogeneous
system~$(\mathcal{ E}_0)$ in which the second members $F_{k_1, \dots,
k_\kappa}^j$ vanish identically. Its solutions are the graphs of
polynomial maps of degree $\leq (\kappa-1)$ from $\K^n$ to $\K^m$. The
defining equations of its Lie algebra of infinitesimal symmetries are
${\bf R}_{ k_1, \dots, k_\kappa }^j =0$, after having replaced the
variables $U_{ l_1, \dots,l_{ \kappa }}^{ i_1}$ by $0= F_{l_1, \dots,
l_\kappa}^{i_1}$ in $I_8$ and $I_9$ in~(\ref{e217}).  We will keep in
this system the only equations coming from the vanishing of the
coefficients of the five families of monomials ${\rm ct.}$,
$U_{l_1,\dots,l_{\kappa-2 }}^{i_1}$, $U_{l_1, \dots, l_{\kappa-1
}}^{i_1}$, $U_{l_1 }^{i_1}\, U_{l_2, \dots, l_{\kappa}}^{i_2}$ and
$U_{l_1,l_2}^{i_1}\, U_{l_3, \dots, l_{\kappa+1 }}^{i_2}$ (this is
inspired from the computations in Subsection~\ref{sub3.1}). The
coefficients of these five monomials families already appear in the
expression~(\ref{e217}). Moreover we fix $l_1=l_2= \cdots=l_{
\kappa+1}=l$ and $i_1=i_2$, except for the coefficient of the monomial
$U_l^{i_1}\, U_{l, \dots, l}^{i_2}$, where we fix first $i_1=i_2$ and
then $i_1 \neq i_2$. This provides the six partial differential linear
equations:
\begin{equation}\label{ee21}
\left\{
\aligned
0 = 
& \
R_{x_{k_1}x_{k_2}\cdots x_{k_{\kappa}}}^j, \\
0 =
& \
\sum_{\sigma\in\mathfrak{S}_\kappa^{\kappa-2}}\,
\delta_{k_{\sigma(1)},\dots,k_{\sigma(\kappa-2)}}^{l,\dots\dots\dots,l}\,
\
R_{x_{k_{\sigma(\kappa-1)}}x_{k_{\sigma(\kappa)}}u^{i_1}}^j- \\
& \ \ \ \ \ \ \ \ \ \ \ \ \ \ \ \ \ \ \ \ \ \ \ \ -
\delta_{i_1}^j\left(
\sum_{\sigma\in\mathfrak{S}_\kappa^{\kappa-3}}\, 
\delta_{k_{\sigma(1)},\dots,k_{\sigma(\kappa-3)}}^{l,\dots\dots\dots,l}\, 
Q_{x_{k_{\sigma(\kappa-2)}}x_{k_{\sigma(\kappa-1)}}x_{k_{\sigma(k)}}}^l
\right), \\
0 =
& \
\sum_{\sigma\in\mathfrak{S}_\kappa^{\kappa-1}}\, 
\delta_{k_{\sigma(1)},\dots,k_{\sigma(\kappa-1)}}^{l,\dots\dots\dots,l}\, 
R_{x_{k_{\sigma(\kappa)}}u^{i_1}}^j- \\
& \ \ \ \ \ \ \ \ \ \ \ \ \ \ \ \ \ \ \ \ \ \ \ \ -
\delta_{i_1}^j\left(
\sum_{\sigma\in\mathfrak{S}_\kappa^{\kappa-2}}\, 
\delta_{k_{\sigma(1)},\dots,k_{\sigma(\kappa-2)}}^{l,\dots\dots\dots,l}\, 
Q_{x_{k_{\sigma(\kappa-1)}}x_{k_{\sigma(\kappa)}}}^l
\right), \\
0 =
& \
\kappa\, \delta_{k_1,\dots,k_\kappa}^{l,\dots\dots,l}\, 
R_{u^{i_1}u^{i_1}}^j-\kappa \, \delta_{i_1}^j\left(
\sum_{\sigma\in\mathfrak{S}_\kappa^{\kappa-1}}\, 
\delta_{k_{\sigma(1)},\dots,k_{\sigma(\kappa-1)}}^{l,\dots\dots\dots,l}\,
Q_{x_{k_{\sigma(\kappa)}}u^{i_1}}^l\right), \\
0 =
& \
2\kappa \, \delta_{k_1,\dots,k_\kappa}^{l,\dots\dots,l}\, 
R_{u^{i_1}u^{i_2}}^j-\kappa \, \delta_{i_1}^j\left(
\sum_{\sigma\in\mathfrak{S}_\kappa^{\kappa-1}}\, 
\delta_{k_{\sigma(1)},\dots,k_{\sigma(\kappa-1)}}^{l,\dots\dots\dots,l}\,
Q_{x_{k_{\sigma(\kappa)}}u^{i_2}}^l
\right) - \\
& \ \ \ \ \ \ \ \ \ \ \ \ \ \ \ \ \ \ \ \ \ \ \ \ -
\kappa \, \delta_{i_2}^j \,
\left(
\sum_{\sigma\in\mathfrak{S}_\kappa^{\kappa-1}}\, 
\delta_{k_{\sigma(1)},\dots,k_{\sigma(\kappa-1)}}^{l,\dots\dots\dots,l}\,
Q_{x_{k_{\sigma(\kappa)}}u^{i_1}}^l
\right), \ \ \ \ \ \ \ i_1\neq i_2, \\
0 =
& \
-C_{\kappa+1}^2\, \delta_{i_1}^j\, 
\delta_{k_1,\dots,k_\kappa}^{l,\dots\dots,l}\, Q_{u^{i_1}}^l.
\endaligned\right.
\end{equation}
To solve the system~(\ref{ee21}) we fix the indices
$k_1=\cdots=k_\kappa=l$ and $j=i_1$ in the sixth equation, implying
$Q_{u^{i_1}}^l=0$. Hence the terms following $\delta_{i_1}^j$ and
$\delta_{i_2}^j$ in the fourth and in the fifth equations vanish
identically. Let us choose the indices $k_1=\dots=k_\kappa$ in the
fourth and the fifth equations (this last equation is satisfied only
for $i_1\neq i_2$). We obtain first three simple equations, without
any restriction on the indices:
\begin{equation}\label{ee22}
\left\{
\aligned
0 = 
& \
R_{x_{k_1}x_{k_2}\cdots x_{k_\kappa}}^j, \\
0 =
& \
Q_{u^{i_1}}^l, \\
0 =
& \
R_{u^{i_1}u^{i_2}}^j.
\endaligned\right.
\end{equation}
Finally we specify the indices in the third equation of~(\ref{ee21})
as follows: $l=k_\kappa=\cdots=k_3=k_2=k_1$~; then
$l=k_\kappa=\cdots=k_3=k_2\neq k_1$~; finally $l=k_\kappa=\cdots=k_3$,
$k_3\neq k_2$, $k_3\neq k_1$. This gives the three following
equations:
\begin{equation}\label{ee23}
\left\{
\aligned
0 =
& \
C_\kappa^1 \, R_{x_{k_1}u^{i_1}}^j -C_\kappa^2 \, \delta_{i_1}^j\, 
Q_{x_{k_1}x_{k_1}}^{k_1}, \\
0 =
& \
R_{x_{k_1}u^{i_1}}^ j -C_{\kappa-1}^1 \, \delta_{i_1}^j \, 
Q_{x_{k_1}x_{k_2}}^{k_2}, \ \ \ \ \ \ \ k_2\neq k_1, \\
0 = 
& \
-\delta_{i_1}^j\, Q_{x_{k_1}x_{k_2}}^{k_3}, \ \ \ \ \ \ \ 
k_3\neq k_1, \ \ k_3\neq k_2.
\endaligned\right.
\end{equation}
We specify the indices in the second equation of~(\ref{ee21}) as follows: 
$l=k_\kappa=\cdots=k_3=k_2=k_1$~; then
$l=k_\kappa=\cdots=k_3=k_2\neq k_1$~; then
$l=k_\kappa=\cdots=k_3$, $k_3\neq k_2$, $k_3\neq k_1$~; finally
$l=k_\kappa=\cdots=k_4$, $l\neq k_1$, $l\neq k_2$, $l\neq k_3$. 
This gives the four following equalities:
\begin{equation}\label{ee24}
\left\{
\aligned
0 =
& \
C_\kappa^2 \, R_{x_{k_1}x_{k_1}u^{i_1}}^j -C_\kappa^3 \, \delta_{i_1}^j\, 
Q_{x_{k_1}x_{k_1}x_{k_1}}^{k_1}, \\
0 =
& \
C_{\kappa-1}^1 \, R_{x_{k_1}x_{k_2}u^{i_1}}^j -
C_{\kappa-1}^2 \, \delta_{i_1}^j\, Q_{x_{k_1}x_{k_2}x_{k_2}}^{k_2}, 
\ \ \ \ \ \ \ k_2\neq k_1, \\
0 =
& \
R_{x_{k_1}x_{k_2}u^{i_1}}^j -  C_{\kappa-2}^1\, \delta_{i_1}^j \, 
Q_{x_{k_1}x_{k_2}x_{k_3}}^{k_3}, \ \ \ \ \ \ \ 
k_3\neq k_1, \ \ k_3 \neq k_2, \\
0 =
& \
-\delta_{i_1}^j\, Q_{x_{k_1}x_{k_2}x_{k_3}}^l, \ \ \ \ \ \ \ 
l\neq k_1, \ \ l\neq k_2, \ \ l\neq k_3.
\endaligned\right.
\end{equation}
Let us differentiate now the equations~(\ref{ee23}) with respect to
the variables $x_l$ as follows: we differentiate (\ref{ee23})$_1$ with
respect to $x_{k_1}$~; then we differentiate (\ref{ee23})$_2$ with
respect to $x_{k_2}$~; finally we differentiate~ (\ref{ee23})$_3$ with
respect to $x_{k_3}$. This gives the three following equations:
\begin{equation}\label{ee25}
\left\{
\aligned
0 =
& \
C_\kappa^1 \, R_{x_{k_1}x_{k_1}u^{i_1}}^j-
C_\kappa^2\, \delta_{i_1}^j\, Q_{x_{k_1}x_{k_1}x_{k_1}}^{k_1}, \\
0 = 
& \
R_{x_{k_1}x_{k_2}u^{i_1}}^j-C_{\kappa-1}^1 \, \delta_{i_1}^j\, 
Q_{x_{k_1}x_{k_2}x_{k_2}}^{k_2}, \ \ \ \ \ \ \ 
k_2\neq k_1, \\
0 =
& \
-\delta_{i_1}^j\, Q_{x_{k_1}x_{k_2}x_{k_3}}^{k_3}, \ \ \ \ \ \ \ 
k_3\neq k_1, \ \ k_3\neq k_2.
\endaligned\right.
\end{equation}
The seven equations given by the systems~(\ref{ee24}) and~(\ref{ee25})
may be considered as three systems of two equations (of two variables)
with a nonzero determinant, to which we add the last equation
(\ref{ee24})$_4$. We get immediately:
\begin{equation}\label{ee26}
\left\{
\aligned
0 = 
& \ 
R_{x_{k_1}x_{k_1}u^{i_1}}^j=
\delta_{i_1}^j\, Q_{x_{k_1}x_{k_1}x_{k_1}}^{k_1}, \\
0 =
& \
R_{x_{k_1}x_{k_2}u^{i_1}}^j = 
\delta_{i_1}^j\, Q_{x_{k_1}x_{k_2}x_{k_2}}^{k_2}, 
\ \ \ \ \ \ \ k_2\neq k_1, \\
0 =
& \
R_{x_{k_1}x_{k_2}u^{i_1}}^j =
\delta_{i_1}^j\, Q_{x_{k_1}x_{k_2}x_{k_3}}^{k_3}, 
\ \ \ \ \ \ \ k_3\neq k_1, \ \ k_3\neq k_2, \\
0 =
& \ \delta_{i_1}^j \, Q_{x_{k_1}x_{k_2}x_{k_3}}^l, 
\ \ \ \ \ \ \ l\neq k_1, \ \ l\neq k_2, \ \ l\neq k_3.
\endaligned\right.
\end{equation} 
It follows from these relations and from the relations
$Q_{u^{i_1}}^l=R_{u^{i_1}u^{i_2}}^j=0$ obtained in~(\ref{ee22}) that
all the third order partial derivatives of $Q^l$ vanish identically,
this being also satisfied by the third order partial derivatives of
$R^j$ containing at least one partial derivative with respect to
$u^{i_1}$:
\begin{equation}\label{ee27}
\left\{
\aligned
0 = 
& \ 
Q_{x_{k_1}x_{k_2}x_{k_3}}^l=
Q_{x_{k_1}x_{k_2}u^{i_1}}^l=
Q_{x_{k_1}u^{i_1}u^{i_2}}^l=
Q_{u^{i_1}u^{i_2}u^{i_3}}^l, \\
0 =
& \
R_{x_{k_1}x_{k_2}u^{i_1}}^j=
R_{x_{k_1}u^{i_1}u^{i_2}}^j=
R_{u^{i_1}u^{i_2}u^{i_3}}^j.
\endaligned\right.
\end{equation} 
It follows from the equations~(\ref{ee22}) and~(\ref{ee27}) that all
the functions $Q^l$ are polynomials of degree $\leq 2$ with respect to
the variables $x_{k_1}$ and all the functions $R^j$ are a sum of a
polynomial of degree $\leq (\kappa-1)$ in the variables $x_{k_1}$ and
of monomials of the form $u^{i_1}$ and $x_{k_1}u^{i_1}$.  Let us
develop now the relations~(\ref{ee23}) separately for $j=i_1$ and
$j\neq i_1$. We obtain the five equations:
\begin{equation}\label{ee28}
\left\{
\aligned
0 =
& \
C_\kappa^1\, R_{x_{k_1}u^{i_1}}^{i_1}-C_\kappa^2\, 
Q_{x_{k_1}x_{k_1}}^{k_1}, \\
0 =
& \
C_\kappa^1 \, R_{x_{k_1}u^{i_1}}^j, \ \ \ \ \ \ \ j\neq i_1, \\
0 = 
& \
R_{x_{k_1}u^{i_1}}^{i_1}-C_{\kappa-1}^1\, Q_{x_{k_1}x_{k_2}}^{k_2}, 
\ \ \ \ \ \ \ k_2\neq k_1, \\
0 =
& \
R_{x_{k_1}u^{i_1}}^j, \ \ \ \ \ \ \ j\neq i_1, \\
0 =
& \
-Q_{x_{k_1}x_{k_2}}^{k_3}, \ \ \ \ \ \ \ k_3\neq k_1, \ \ k_3\neq k_2.
\endaligned\right.
\end{equation}
According to the equations~(\ref{ee22}), (\ref{ee27}), (\ref{ee28}), 
we have the following form of the general solution:
\begin{equation}\label{e47}
\left\{
\aligned
Q^l(x,u)=
& \
A^l+\sum_{k_1=1}^n\, B_{k_1}^l\, x_{k_1}+
\sum_{k_1=1}^n\, C_{k_1}\, x_{k_1}\, x_l, \\
R^j(x,u) =
& \
\sum_{k_1=1}^n\, 
(\kappa-1) \, C_{k_1} \, x_{k_1} \, 
u^j +
\sum_{i_1=1}^m\, D_{i_1}^j\, u^{i_1} +
E^{j,0}+
\sum_{k_1=1}^n\, 
E_{k_1}^{j,1}\, x_{k_1}+ \\
& \
+\cdots+
\sum_{1\leq k_1\leq \cdots \leq k_{\kappa-1}\leq n}\, 
E_{k_1,\dots,k_{\kappa-1}}^{j,\kappa-1}\, 
x_{k_1}\cdots x_{k_{\kappa-1}}.
\endaligned\right.
\end{equation}
Here the $n+n^2+n+m^2+m\, C_{n+\kappa-1}^{\kappa-1}$ constants $A^l,
\, B_{k_1}^l, \, C_{k_1}, \, D_{i_1}^j, \, E^{j,0}, \, E_{k_1}^{j,1}$,
$\dots$, $E_{k_1,\dots,k_{\kappa-1}}^{j,\kappa-1}\in\K$ are
arbitrary. Moreover one can check that the vector space spanned by the
vector fields
\begin{equation}\label{e48}
\left\{
\aligned
{}
&
\frac{\partial }{\partial x_{k_1}}, \ \ 
x_{k_1}\,\frac{\partial }{\partial x_{k_2}}, \\
&
x_{k_1}\left(
x_1\, \frac{\partial }{\partial x_1}+\cdots+
x_n\, \frac{\partial }{\partial x_n}+
(\kappa-1)\left(
u^1\, \frac{\partial }{\partial u^1}+\cdots+
u^m\, \frac{\partial }{\partial u^m}
\right)
\right),\\
&
u^{i_1}\, \frac{\partial }{\partial u^{i_2}}, \ \ 
\frac{\partial }{\partial u^{i_1}}, \ \
x_{k_1}\, \frac{\partial }{\partial u^{i_1}}, \ \dots\dots, \
x_{k_1}\cdots x_{k_{\kappa-1}}\,  
\frac{\partial }{\partial u^{i_1}},
\endaligned\right.
\end{equation}
is stable under the Lie bracket action and that the flow of each of
these generators is indeed a Lie symmetry of the system
$(\mathcal{E}_0)$. Finally the Lie symmetries of $(\mathcal E_0)$ have
the following form:
\begin{equation}\label{e49}
{\footnotesize
\aligned
& (x,u)\longmapsto
\left(
\left(
\frac{\alpha_{l,0}+\sum_{k=1}^n\, \alpha_{l,k}\, 
x_k}{1+\sum_{k=1}^n\, \varepsilon_k \, x_k}
\right)_{1\leq l\leq n}, \right.\\
& \left.
\left(
\frac{
\sum_{i_1=1}^m\, 
\beta_{i_1}^j\, u_{i_1}+
\gamma^{0,j}+
\sum_{k_1=1}^n\, \gamma_{k_1}^{1,j}\, x_{k_1}+\cdots+
\sum_{k_1\leq \cdots \leq k_{\kappa-1}}\, 
\gamma_{k_1,\dots,k_{\kappa-1}}^{\kappa-1,j}\, 
x_{k_1}\cdots x_{k_{\kappa-1}}}{
[1+\sum_{k=1}^n\, \varepsilon_k \, x_k]^{\kappa-1}}
\right)_{1\leq j\leq m} 
\right). 
\endaligned
}
\end{equation}

\noindent We note again that these transformations preserve the
solutions of $(\mathcal{E}_0): \ \ u_{x_{k_1}\cdots
x_{k_\kappa}}^j=0$, namely the graphs of polynomial maps of degree
$\leq (\kappa-1)$ from $\K^n$ to $\K^m$.

\subsection{Nonhomogeneous system} Let $\kappa \geq 3$. Let us expand 
the defining equations~(\ref{e41}) as done in~(\ref{e42}). We will 
write only the coefficients of the five monomial families ${\rm ct.}$, 
$U_{ l_1, \dots, l_{
\kappa-2}}^{ i_1}$, $U_{ l_1, \dots, l_{ \kappa-1}}^{ i_1}$, $U_{ l_1
}^{ i_1}\, U_{ l_2, \dots, l_{ \kappa}}^{i_2}$ and
$U_{l_1,l_2}^{i_1}\, U_{ l_3, \dots, l_{ \kappa+1 }}^{
i_2}$. Moreover, we fix always $l_1= l_2= \dots = l_\kappa= l_{
\kappa+1} =l$ and $i_1 =i_2$, except for the fourth family of
monomials where we distinguish the two cases $i_1= i_2$ and $i_1 \neq
i_2$. Thus we obtain six linear equations of partial derivatives, the
members on the left side (coming from the expression of ${\bf R}_{
k_1, \dots, k_\kappa}^{j_1}$ given by Lemma~\ref{lem22}) coincide with
the members on the right hand side of~(\ref{ee21}).  Furthermore, the
members on the right hand side are exactly the same as those obtained
in~(\ref{e310}), with more indices! We use the letters
$l',k_1',\dots,k_\kappa'=1,\dots,n$ and $j',i_1'=1,\dots,m$ for the
indices of the arguments of the expressions $\Pi$, obtaining the six
following equations, which generalize the equations~(\ref{e310}):

\begin{equation}\label{e413}
[1]: \ \ \ \ \
R_{x_{k_1}x_{k_2}\cdots x_{k_{\kappa}}}^j= 
\Pi\left(x,u,Q^{l'},Q_{x_{k_1'}}^{l'}, R^{j'}, R_{x_{k_1'}}^{j'}, \dots,
R_{x_{k_1'}\cdots x_{k_{\kappa-1}'}}^{j'}, R_{u^{i_1'}}^{j'}\right),
\end{equation}

$$
\aligned
{}
&
[1]: \ \ \ \ \
R_{x_{k_1}x_{k_2}\cdots x_{k_{\kappa}}}^j= 
\Pi\left(x,u,Q^{l'},Q_{x_{k_1'}}^{l'}, R^{j'}, R_{x_{k_1'}}^{j'}, \dots,
R_{x_{k_1'}\cdots x_{k_{\kappa-1}'}}^{j'}, R_{u^{i_1'}}^{j'}\right).
\\
&
[2]: \ \ \ \ \ 
\sum_{\sigma\in\mathfrak{S}_\kappa^{\kappa-2}}\,
\delta_{k_{\sigma(1)},\dots,k_{\sigma(\kappa-2)}}^{l,\dots\dots\dots,l}\,
R_{x_{k_{\sigma(\kappa-1)}}x_{k_{\sigma(\kappa)}}u^{i_1}}^j - \\
& \ \ \ \ \ \ \ \ \ \ \ \ \ \ \ \ \ \ \ \ \ \ \ \ \ -
\delta_{i_1}^j\left(
\sum_{\sigma\in\mathfrak{S}_\kappa^{\kappa-3}}\, 
\delta_{k_{\sigma(1)},\dots,k_{\sigma(\kappa-3)}}^{l,\dots\dots\dots,l}\, 
Q_{x_{k_{\sigma(\kappa-2)}}x_{k_{\sigma(\kappa-1)}}x_{k_{\sigma(k)}}}^l
\right)= \\
&
= 
\Pi\left(x,u,Q^{l'},Q_{x_{k_1'}}^{l'},Q_{x_{k_1'}x_{k_2'}}^{l'}, 
R^{j'}, R_{x_{k_1'}}^{j'}, \dots,
R_{x_{k_1'}\cdots x_{k_{\kappa-1}'}}^{j'}, R_{u^{i_1'}}^{j'}, 
R_{x_{k_1'}u^{i_1'}}^{j'}\right).
\endaligned
$$

$$
\aligned
{}
&
[3]: \ \ \ \ \ 
\sum_{\sigma\in\mathfrak{S}_\kappa^{\kappa-1}}\, 
\delta_{k_{\sigma(1)},\dots,k_{\sigma(\kappa-1)}}^{l,\dots\dots\dots,l}\, 
R_{x_{k_{\sigma(\kappa)}}u^{i_1}}^j- \\
& \ \ \ \ \ \ \ \ \ \ \ \ \ \ \ \ \ \ \ \ \ \ \ \ \ \ \ \ \ \ \ \ \ \ \ \ \ \ -
\delta_{i_1}^j\left(
\sum_{\sigma\in\mathfrak{S}_\kappa^{\kappa-2}}\, 
\delta_{k_{\sigma(1)},\dots,k_{\sigma(\kappa-2)}}^{l,\dots\dots\dots,l}\, 
Q_{x_{k_{\sigma(\kappa-1)}}x_{k_{\sigma(\kappa)}}}^l
\right)= \\
&
=
\Pi
\left(x,u,Q^{l'},Q_{x_{k_1'}}^{l'},
R^{j'}, R_{x_{k_1'}}^{j'}, \dots,
R_{x_{k_1'}\cdots x_{k_{\kappa-1}'}}^{j'}, R_{u^{i_1'}}^{j'}
\right). \\
\endaligned
$$

$$
\aligned
{}
&
[4]: \ \ \ \ \ 
\kappa\, \delta_{k_1,\dots,k_\kappa}^{l,\dots\dots,l}\, 
R_{u^{i_1}u^{i_1}}^j-\kappa \, \delta_{i_1}^j\left(
\sum_{\sigma\in\mathfrak{S}_\kappa^{\kappa-1}}\, 
\delta_{k_{\sigma(1)},\dots,k_{\sigma(\kappa-1)}}^{l,\dots\dots\dots,l}\,
Q_{x_{k_{\sigma(\kappa)}}u^{i_1}}^l\right)= \\
&
=
\Pi
\left(
x,u,Q^{l'},Q_{x_{k_1'}}^{l'},\dots,Q_{x_{k_1'}\cdots x_{k_{\kappa-1}'}}^{l'},
Q_{u^{i_1'}}^{l'},
R^{j'}, R_{x_{k_1'}}^{j'}, \dots,
R_{x_{k_1'}\cdots x_{k_{\kappa-1}'}}^{j'}, \right. \\
& \ \ \ \ \ \ \ \ \ \ \ \ \ \ \ \ \ \ \ \ \ \ \ \ \ \ \ \ \ \ 
\left.
, R_{u^{i_1'}}^{j'}, 
R_{x_{k_1'}u^{i_1'}}^{j'},\dots,
R_{x_{k_1'}\cdots x_{k_{\kappa-2}'}u^{i_1'}}^{j'}
\right). \\
\endaligned
$$

$$
\aligned
{}
&
[5]: \ \ \ \ \ 
2\kappa \, \delta_{k_1,\dots,k_\kappa}^{l,\dots\dots,l}\, 
R_{u^{i_1}u^{i_2}}^j-\kappa \, \delta_{i_1}^j\left(
\sum_{\sigma\in\mathfrak{S}_\kappa^{\kappa-1}}\, 
\delta_{k_{\sigma(1)},\dots,k_{\sigma(\kappa-1)}}^{l,\dots\dots\dots,l}\,
Q_{x_{k_{\sigma(\kappa)}}u^{i_2}}^l
\right) - \\
& \ \ \ \ \ \ \ \ \ \ \ \ \ \ \ \ \ \ \ \ \ \ \ \ \ \ \ \ -
\kappa \, \delta_{i_2}^j \,
\left(
\sum_{\sigma\in\mathfrak{S}_\kappa^{\kappa-1}}\, 
\delta_{k_{\sigma(1)},\dots,k_{\sigma(\kappa-1)}}^{l,\dots\dots\dots,l}\,
Q_{x_{k_{\sigma(\kappa)}}u^{i_1}}^l
\right), \ \ \ \ \ \ \ i_1\neq i_2. \\
&
=
\Pi
\left(
x,u,Q^{l'},Q_{x_{k_1'}}^{l'},\dots,Q_{x_{k_1'}\cdots x_{k_{\kappa-1}'}}^{l'},
Q_{u^{i_1'}}^{l'},
R^{j'}, R_{x_{k_1'}}^{j'}, \dots,
R_{x_{k_1'}\cdots x_{k_{\kappa-1}'}}^{j'}, \right. \\
& \ \ \ \ \ \ \ \ \ \ \ \ \ \ \ \ \ \ \ \ \ \ \ \ \ \ \ \ \ \ 
\left.
, R_{u^{i_1'}}^{j'}, 
R_{x_{k_1'}u^{i_1'}}^{j'},\dots,
R_{x_{k_1'}\cdots x_{k_{\kappa-2}'}u^{i_1'}}^{j'}
\right). \\
\endaligned
$$

$$
\aligned
{}
&
[6]: \ \ \ \ \ 
-C_{\kappa+1}^2\, \delta_{i_1}^j\, 
\delta_{k_1,\dots,k_\kappa}^{l,\dots\dots,l}\, Q_{u^{i_1}}^l= \\
& 
=
\Pi
\left(
x,u,Q^{l'},Q_{x_{k_1'}}^{l'},\dots,Q_{x_{k_1'}\cdots x_{k_{\kappa-2}'}}^{l'},
R^{j'}, R_{x_{k_1'}}^{j'}, \dots,
R_{x_{k_1'}\cdots x_{k_{\kappa-1}'}}^{j'}, \right. \\
& \ \ \ \ \ \ \ \ \ \ \ \ \ \ \ \ \ \ \ \ \ \ \ \ \ \ \ \ \ \ 
\left.
, R_{u^{i_1'}}^{j'}, 
R_{x_{k_1'}u^{i_1'}}^{j'},\dots,
R_{x_{k_1'}\cdots x_{k_{\kappa-3}'}u^{i_1'}}^{j'}
\right).
\endaligned
$$
\noindent
Then we get the following Lemma:

\begin{lemma}\label{lem41}
Let $J$ denote the collection of $n+n^2+n+m \, C_{n+\kappa-1}^{\kappa-1}+ m^2$
partial derivatives
\begin{equation}\label{e414}
J:= \left(Q^{l'},Q_{x_{k_1'}}^{l'},Q_{x_{k_1'}x_{k_1'}}^{k_1'},
R^{j'},R_{x_{k_1'}}^{j'},\dots,
R_{x_{k_1'}\cdots x_{k_{\kappa-1}'}}^{j'},R_{u^{i_1'}}^{j'}\right).
\end{equation}
After linear combinations on the system~(\ref{e413}) we obtain the 
following equations:
\begin{equation}\label{e415}
\left\{
\aligned
\Pi(x,u,J) =
& \
R_{x_{k_1}\cdots x_{k_\kappa}}^j, \\
\Pi(x,u,J) =
& \ 
Q_{u^{i_1}}^l,
\\
\Pi(x,u,J) =
& \ 
R_{u^{i_1}u^{i_2}}^j,
\\
\Pi(x,u,J) =
& \ 
Q_{x_{k_1}x_{k_2}x_{k_3}}^l,
\\
\Pi(x,u,J) =
& \ 
R_{x_{k_1}u^{i_1}}^j,
\\
\Pi(x,u,J) =
& \ 
Q_{x_{k_1}x_{k_2}}^{k_1}, \ \ \ \ \ \ \ k_1\neq k_2,
\\
\Pi(x,u,J) =
& \ 
Q_{x_{k_1}x_{k_2}}^l, \ \ \ \ \ \ \ l\neq k_1, \ \ l\neq k_2.
\endaligned\right.
\end{equation}
Moreover all the partial derivatives (with respect to $x_l$ and $u^i$)
up to order three of the coefficients $Q^l$ and $R^j$ of the vector
field $X\in \mathfrak{ Sym} (\mathcal{ E})$ are of the form
$\Pi(x,u,J)$. Hence every function $Q^l$ and $R^j$ is uniquely
determined by the values at the origin of the $n+ n^2+n+ m \, C_{n+
\kappa -1}^{ \kappa-1}+ m^2$ partial derivatives~(\ref{e414}).  This
implies that $\dim_\K \, \mathfrak{ Sym} (\mathcal{E}) \leq n^2+ 2n+
m^2+ m\, C_{ n+ \kappa-1}^{ \kappa-1}$.
\end{lemma}

\noindent{\it Proof}. Since the second part of Lemma~\ref{lem41} is
immediate let us establish only the identities~(\ref{e415}). We first
specify the indices in the equation~(\ref{e413})$_{[3]}$ as follows:
$l=k_\kappa=\cdots=k_3=k_2=k_1$~; then
$l=k_\kappa=\cdots=k_3=k_2\neq k_1$~; and finally
$l=k_\kappa=\cdots=k_3$, $k_3\neq k_2$, $k_3\neq k_1$. This gives 
three equations whose members on the right hand side are the same 
as those in the equation~(\ref{ee23}) and whose members on the left 
hand side are the same as those in the equation~(\ref{e413})$_{[3]}$:
\begin{equation}\label{e416}
{\small 
\left\{
\aligned
{}
&
\Pi
\left(x,u,Q^{l'},Q_{x_{k_1'}}^{l'},
R^{j'}, R_{x_{k_1'}}^{j'}, \dots,
R_{x_{k_1'}\cdots x_{k_{\kappa-1}'}}^{j'}, R_{u^{i_1'}}^{j'}
\right)=
C_\kappa^1 \, R_{x_{k_1}u^{i_1}}^j -C_\kappa^2 \, \delta_{i_1}^j\, 
Q_{x_{k_1}x_{k_1}}^{k_1},
\\
&
\Pi
\left(x,u,Q^{l'},Q_{x_{k_1'}}^{l'},
R^{j'}, R_{x_{k_1'}}^{j'}, \dots,
R_{x_{k_1'}\cdots x_{k_{\kappa-1}'}}^{j'}, R_{u^{i_1'}}^{j'}
\right)=
R_{x_{k_1}u^{i_1}}^ j -C_{\kappa-1}^1 \, \delta_{i_1}^j \, 
Q_{x_{k_1}x_{k_2}}^{k_2}, 
\ \ k_2\neq k_1, \\
&
\Pi
\left(x,u,Q^{l'},Q_{x_{k_1'}}^{l'},
R^{j'}, R_{x_{k_1'}}^{j'}, \dots,
R_{x_{k_1'}\cdots x_{k_{\kappa-1}'}}^{j'}, R_{u^{i_1'}}^{j'}
\right)=
-\delta_{i_1}^j\, Q_{x_{k_1}x_{k_2}}^{k_3}, 
\ \ \ \ \ \ \ k_3\neq k_1, \ \ k_3\neq k_2.
\endaligned\right.
}
\end{equation}

\noindent
We remark that these three equations (after specialization of $j=i_1$
or of $j\neq i_1$ and after some easy linear combinations) provide
directly the fifth, sixth and seventh equations of~(\ref{e415}).  In
particular we may replace the values of the partial derivatives
$R_{x_{j_1'}u^{i_1'}}^{j'}$ and $Q_{x_{k_1'}x_{k_2'}}^{l'}$ with
$k_1'\neq k_2'$ or $l'\neq k_1',\, l'\neq k_2'$ appearing in the
expressions $\Pi$ of the second member of~(\ref{e413})$_{[1]}$ by
their values just obtained from the fifth, the sixth and the seventh
equations of~(\ref{e415}). This gives the first equation
of~(\ref{e415}).

Then we specify the indices in (\ref{e413})$_{[2]}$ as follows:
$l=k_\kappa=\cdots=k_3=k_2=k_1$~; then $l=k_\kappa=\cdots=k_3=k_2\neq
k_1$~; then $l=k_\kappa=\cdots=k_3$, $k_3\neq k_2$, $k_3\neq k_1$~;
and finally $l=k_\kappa=\cdots=k_4$, $l\neq k_1$, $l\neq k_2$, $l\neq
k_3$. This gives four equations, whose members on the right hand side
are the same as those in~(\ref{ee24}) and the members on the left hand
side are the same as those in~(\ref{e413})$_{[2]}$:

\begin{equation}\label{e417}
{\small
\left\{
\aligned
{}
&
\Pi\left(x,u,Q^{l'},Q_{x_{k_1'}}^{l'},Q_{x_{k_1'}x_{k_2'}}^{l'}, 
R^{j'}, R_{x_{k_1'}}^{j'}, \dots,
R_{x_{k_1'}\cdots x_{k_{\kappa-1}'}}^{j'}, R_{u^{i_1'}}^{j'}, 
R_{x_{k_1'}u^{i_1'}}^{j'}\right)= \\
& \ \ \ \ \ \ \ \ \ \ \ \ \ \ \ \
=C_\kappa^2 \, R_{x_{k_1}x_{k_1}u^{i_1}}^j -C_\kappa^3 \, \delta_{i_1}^j\, 
Q_{x_{k_1}x_{k_1}x_{k_1}}^{k_1}, \\
&
\Pi\left(x,u,Q^{l'},Q_{x_{k_1'}}^{l'},Q_{x_{k_1'}x_{k_2'}}^{l'}, 
R^{j'}, R_{x_{k_1'}}^{j'}, \dots,
R_{x_{k_1'}\cdots x_{k_{\kappa-1}'}}^{j'}, R_{u^{i_1'}}^{j'}, 
R_{x_{k_1'}u^{i_1'}}^{j'}\right)= \\
& \ \ \ \ \ \ \ \ \ \ \ \ \ \ \ \
=C_{\kappa-1}^1 \, R_{x_{k_1}x_{k_2}u^{i_1}}^j -
C_{\kappa-1}^2 \, \delta_{i_1}^j\, Q_{x_{k_1}x_{k_2}x_{k_2}}^{k_2}, 
\ \ \ \ \ \ \ k_2\neq k_1, \\
&
\Pi\left(x,u,Q^{l'},Q_{x_{k_1'}}^{l'},Q_{x_{k_1'}x_{k_2'}}^{l'}, 
R^{j'}, R_{x_{k_1'}}^{j'}, \dots,
R_{x_{k_1'}\cdots x_{k_{\kappa-1}'}}^{j'}, R_{u^{i_1'}}^{j'}, 
R_{x_{k_1'}u^{i_1'}}^{j'}\right)= \\
& \ \ \ \ \ \ \ \ \ \ \ \ \ \ \ \
=R_{x_{k_1}x_{k_2}u^{i_1}}^j -  C_{\kappa-2}^1\, \delta_{i_1}^j \, 
Q_{x_{k_1}x_{k_2}x_{k_3}}^{k_3}, \ \ \ \ \ \ \ 
k_3\neq k_1, \ \ k_3 \neq k_2, \\
&
\Pi\left(x,u,Q^{l'},Q_{x_{k_1'}}^{l'},Q_{x_{k_1'}x_{k_2'}}^{l'}, 
R^{j'}, R_{x_{k_1'}}^{j'}, \dots,
R_{x_{k_1'}\cdots x_{k_{\kappa-1}'}}^{j'}, R_{u^{i_1'}}^{j'}, 
R_{x_{k_1'}u^{i_1'}}^{j'}\right)= \\
& \ \ \ \ \ \ \ \ \ \ \ \ \ \ \ \
=-\delta_{i_1}^j\, Q_{x_{k_1}x_{k_2}x_{k_3}}^l, \ \ \ \ \ \ \ 
l\neq k_1, \ \ l\neq k_2, \ \ l\neq k_3.
\endaligned\right.
}
\end{equation}

\noindent
Using the fifth, the sixth and the seventh equations of~(\ref{e415})
just obtained, we may replace the partial derivatives
$R_{x_{j_1'}u^{i_1'}}^{j'}$ and $Q_{x_{k_1'}x_{k_2'}}^{l'}$ with
$k_1'\neq k_2'$ or $l'\neq k_1'$, $l'\neq k_2'$ appearing in the
expressions $\Pi$ of~(\ref{e417}), providing four new equations in
which the arguments of $\Pi$ are the desired ones: $(x,u,J)$, where
$J$ is defined in~(\ref{e414}):
\begin{equation}\label{e418}
\left\{
\aligned
\Pi(x,u,J) =
& \
C_\kappa^2 \, R_{x_{k_1}x_{k_1}u^{i_1}}^j -C_\kappa^3 \, \delta_{i_1}^j\, 
Q_{x_{k_1}x_{k_1}x_{k_1}}^{k_1}, \\
\Pi(x,u,J) =
& \
C_{\kappa-1}^1 \, R_{x_{k_1}x_{k_2}u^{i_1}}^j -
C_{\kappa-1}^2 \, \delta_{i_1}^j\, Q_{x_{k_1}x_{k_2}x_{k_2}}^{k_2}, 
\ \ \ \ \ \ \ k_2\neq k_1, \\
\Pi(x,u,J) =
& \
R_{x_{k_1}x_{k_2}u^{i_1}}^j -  C_{\kappa-2}^1\, \delta_{i_1}^j \, 
Q_{x_{k_1}x_{k_2}x_{k_3}}^{k_3}, \ \ \ \ \ \ \ 
k_3\neq k_1, \ \ k_3 \neq k_2, \\
\Pi(x,u,J) =
& \
-\delta_{i_1}^j\, Q_{x_{k_1}x_{k_2}x_{k_3}}^l, \ \ \ \ \ \ \ 
l\neq k_1, \ \ l\neq k_2, \ \ l\neq k_3.
\endaligned\right.
\end{equation}
\noindent
Let us differentiate now the equations~(\ref{e416}) with respect to
the variables $x_l$ as follows: first we differentiate
(\ref{e416})$_1$ with respect to $x_{k_1}$~; then we differentiate
(\ref{e416})$_2$ with respect to $x_{k_2}$~; finally we differentiate
(\ref{e416})$_3$ with respect to $x_{k_3}$. The arguments in the
expressions $\Pi$ in the equation~(\ref{e416}) contain now the terms
$R_{x_{k_1'}\cdots x_{k_\kappa'}}^{j'}$; we replace them by their
value given in the first equation of~(\ref{e415}) already
obtained. The arguments also contain the terms
$R_{x_{j_1'}u^{i_1'}}^{j'}$ and $Q_{x_{k_1'}x_{k_2'}}^{l'}$ with
$k_1'\neq k_2'$ or $l'\neq k_1'$, $l'\neq k_2'$. We replace them by
their value given by the fifth, the sixth and the seventh equations
of~(\ref{e415}). We obtain three new equations in which the arguments
of the expressions $\Pi$ are the desired ones: $(x,u,J)$, where $J$ is
defined in~(\ref{e414}):
\begin{equation}\label{e419}
\left\{
\aligned
{}
&
\Pi
\left(x,u,J\right)
=C_\kappa^1 \, R_{x_{k_1}x_{k_1}u^{i_1}}^j-
C_\kappa^2\, \delta_{i_1}^j\, Q_{x_{k_1}x_{k_1}x_{k_1}}^{k_1}, \\
&
\Pi
\left(x,u,J\right)=
R_{x_{k_1}x_{k_2}u^{i_1}}^j-C_{\kappa-1}^1 \, \delta_{i_1}^j\, 
Q_{x_{k_1}x_{k_2}x_{k_2}}^{k_2}, \ \ \ \ \ \ \ 
k_2\neq k_1, \\
&
\Pi\left(x,u,J\right)=
-\delta_{i_1}^j\, Q_{x_{k_1}x_{k_2}x_{k_3}}^{k_3}, \ \ \ \ \ \ \ 
k_3\neq k_1, \ \ k_3\neq k_2.
\endaligned\right.
\end{equation}
The seven equations~(\ref{e418}) and~(\ref{e419}) may be considered as
three systems of two linear equations of two variables with a nonzero
determinant, the seventh equation being the last equation
in~(\ref{e418}). We immediately obtain:
\begin{equation}\label{e420}
\left\{
\aligned
\Pi(x,u,J) = 
& \ 
R_{x_{k_1}x_{k_1}u^{i_1}}^j=
\delta_{i_1}^j\, Q_{x_{k_1}x_{k_1}x_{k_1}}^{k_1}, \\
\Pi(x,u,J) =
& \
R_{x_{k_1}x_{k_2}u^{i_1}}^j = 
\delta_{i_1}^j\, Q_{x_{k_1}x_{k_2}x_{k_2}}^{k_2}, 
\ \ \ \ \ \ \ k_2\neq k_1, \\
\Pi(x,u,J) =
& \
R_{x_{k_1}x_{k_2}u^{i_1}}^j =
\delta_{i_1}^j\, Q_{x_{k_1}x_{k_2}x_{k_3}}^{k_3}, 
\ \ \ \ \ \ \ k_3\neq k_1, \ \ k_3\neq k_2, \\
\Pi(x,u,J) =
& \ \delta_{i_1}^j \, Q_{x_{k_1}x_{k_2}x_{k_3}}^l, 
\ \ \ \ \ \ \ k_3\neq k_1, \ \ k_3 \neq k_2,
\endaligned\right.
\end{equation} 
giving the fourth equation in~(\ref{e415}).

It remains now to obtain the second and the third equations
in~(\ref{e415}). Let us write firstly equation~(\ref{e413})$_{[6]}$
with the choice of the indices $j=i_1$, $l=k_1=\dots=k_\kappa$. This
gives the equation:
\begin{equation}\label{e421}
\left\{
\aligned
Q_{u^{i_1}}^l=\Pi
&
\left(
x,u,Q^{l'},
Q_{x_{k_1'}}^{l'},\dots,Q_{x_{k_1'}\cdots x_{k_{\kappa-2}'}}^{l'},
R^{j'}, R_{x_{k_1'}}^{j'}, \dots,
R_{x_{k_1'}\cdots x_{k_{\kappa-1}'}}^{j'}, \right. \\
& \ \ \ \ \ \ \ \ \ \ \ \ \ \ \ \ \ \ \ \ \ \ \ \ \ \ \ \ \ \ 
\left.
, R_{u^{i_1'}}^{j'}, 
R_{x_{k_1'}u^{i_1'}}^{j'},\dots,
R_{x_{k_1'}\cdots x_{k_{\kappa-3}'}u^{i_1'}}^{j'}
\right).
\endaligned\right.
\end{equation}
We observe first that the differentiation with respect to the
variables $x_l$ of one of the expressions $\Pi(x,u,J)$ remains an
expression $\Pi(x,u,J)$. Indeed we see from~(\ref{e414}) that there
appears, in the partial derivative $J_{x_l}$, derivatives
$Q_{x_{k_1'}x_{k_2'}}^{l'}$ with $k_1'\neq k_2'$ or $l'\neq k_1',\,
l'\neq k_2'$. We may replace them by their value obtained in the sixth
and the seventh equations of~(\ref{e415}). It also appears some
derivatives $Q_{x_{k_1'}x_{k_2'}x_{k_3'}}^{l'}$ (we replace them by
their value obtained in the fourth equation of~(\ref{e415})), some
derivatives $R_{x_{k_1'}\cdots x_{k_\kappa'}}^{j'}$ (we replace them
by their value obtained in the first equation of~(\ref{e415})) and
some derivatives $R_{x_{k_1'}u^{i_1'}}^{j'}$ (we replace them by their
value obtained in the fifth equation of~(\ref{e415})). Consequently we
may write:
\begin{equation}\label{e422}
[\Pi(x,u,J)]_{x_l}=\Pi(x,u,J).
\end{equation}
It follows that any derivative with respect to $x_l$ (to any order) of
the fourth and the fifth equations of~(\ref{e415}) provides
expressions of the form $\Pi(x,u,J)$. In other words for any integer
$\lambda\geq 3$ and any integer $\mu\geq 1$ we have
\begin{equation}\label{e423}
\left\{
\aligned
\Pi(x,u,J)=
& \
Q_{x_{k_1}x_{k_2}x_{k_3}\cdots x_{k_\lambda}}^l,\\
\Pi(x,u,J)=
& \
R_{x_{k_1} \cdots x_{k_\mu}u^{i_1}}^j.
\endaligned\right.
\end{equation}

\noindent
We may replace then these values in the equation~(\ref{e421}),
replacing also the derivatives $Q_{x_{k_1'}x_{k_2'}}^{l'}$ with
$k_1'\neq k_2'$ or $l'\neq k_1',\, l'\neq k_2'$ by their values
obtained in the sixth and the seventh equations of~(\ref{e415}). This
gives the second equation of~(\ref{e415}).

We also remark that by a differentiation with respect to the variables
$x_l$, the second equation $Q_{u^{i_1}}^l=\Pi(x,u,J)$ just obtained
implies, using~(\ref{e422}):
\begin{equation}\label{e424}
\Pi(x,u,J)=Q_{x_{k_1}u^{i_1}}^l.
\end{equation}

It remains finally to write (\ref{e413})$_{[4]}$ first with the choice
of indices $l=k_1=\cdots=k_\kappa$, $j=i_1$ then with the choice of
indices $l=k_1=\cdots=k_\kappa$, $j\neq i_1$. We also write
(\ref{e413})$_{[5]}$ first with the choice of indices
$l=k_1=\cdots=k_\kappa$, $j=i_2$ then with the choice of indices
$l=k_1=\cdots=k_\kappa$, $j\neq i_1$, $j\neq i_2$. We obtain four new
equations:
 
\begin{equation}\label{e425}
{\small
\left\{
\aligned
R_{u^{i_1}u^{i_1}}^{i_1}-\kappa \, Q_{x_{k_1}u^{i_1}}^{k_1}
=\Pi
&
\left(
x,u,Q^{l'},
Q_{x_{k_1'}}^{l'},\dots,Q_{x_{k_1'}\cdots x_{k_{\kappa-1}'}}^{l'},
Q_{u^{i_1'}}^{l'},
R^{j'}, R_{x_{k_1'}}^{j'}, \dots,
R_{x_{k_1'}\cdots x_{k_{\kappa-1}'}}^{j'}, \right. \\
& \ \ \ \ \ \ \ \ \ \ \ \ \ \ \ \ \ \ \ \ \ \ \ \ \ \ \ \ \ \ 
\left.
, R_{u^{i_1'}}^{j'}, 
R_{x_{k_1'}u^{i_1'}}^{j'},\dots,
R_{x_{k_1'}\cdots x_{k_{\kappa-2}'}u^{i_1'}}^{j'}
\right), \\
R_{u^{i_1}u^{i_1}}^j
=\Pi
&
\left(
x,u,Q^{l'},
Q_{x_{k_1'}}^{l'},\dots,Q_{x_{k_1'}\cdots x_{k_{\kappa-1}'}}^{l'},
Q_{u^{i_1'}}^{l'},
R^{j'}, R_{x_{k_1'}}^{j'}, \dots,
R_{x_{k_1'}\cdots x_{k_{\kappa-1}'}}^{j'}, \right. \\
& \ \ \ \ \ \ \ \ \ \ \ \ \ \ \ \ \ \ \ \ \ \ \ \ \ \ \ \ \ \ 
\left.
, R_{u^{i_1'}}^{j'}, 
R_{x_{k_1'}u^{i_1'}}^{j'},\dots,
R_{x_{k_1'}\cdots x_{k_{\kappa-2}'}u^{i_1'}}^{j'}
\right), \ \ \ \ \ \ \ j\neq i_1, \\
2R_{u^{i_1}u^{i_2}}^{i_2}-\kappa \, Q_{x_{k_1}u^{i_1}}^{k_1}
=\Pi
&
\left(
x,u,Q^{l'},
Q_{x_{k_1'}}^{l'},\dots,Q_{x_{k_1'}\cdots x_{k_{\kappa-1}'}}^{l'},
Q_{u^{i_1'}}^{l'},
R^{j'}, R_{x_{k_1'}}^{j'}, \dots,
R_{x_{k_1'}\cdots x_{k_{\kappa-1}'}}^{j'}, \right. \\
& \ \ \ \ \ \ \ \ \ \ \ \ \ \ \ \ \ \ \ \ \ \ \ \ \ \ \ \ \ \ 
\left.
, R_{u^{i_1'}}^{j'}, 
R_{x_{k_1'}u^{i_1'}}^{j'},\dots,
R_{x_{k_1'}\cdots x_{k_{\kappa-2}'}u^{i_1'}}^{j'}
\right), \ \ \ \ \ \ \ i_1\neq i_2, \\
R_{u^{i_1}u^{i_2}}^j
=\Pi
&
\left(
x,u,Q^{l'},
Q_{x_{k_1'}}^{l'},\dots,Q_{x_{k_1'}\cdots x_{k_{\kappa-1}'}}^{l'},
Q_{u^{i_1'}}^{l'},
R^{j'}, R_{x_{k_1'}}^{j'}, \dots,
R_{x_{k_1'}\cdots x_{k_{\kappa-1}'}}^{j'}, \right. \\
& \ \ \ 
\left.
, R_{u^{i_1'}}^{j'}, 
R_{x_{k_1'}u^{i_1'}}^{j'},\dots,
R_{x_{k_1'}\cdots x_{k_{\kappa-2}'}u^{i_1'}}^{j'}
\right), \ \ \ \ \ \ \ i_1\neq i_2, \ \ j\neq i_2, \ \ j\neq i_2. \\
\endaligned\right.
}
\end{equation}

\noindent
Using the equations of~(\ref{e415}) we already obtained
(namely all except the second equation), using~(\ref{e423}) 
and~(\ref{e424}), we may simplify these four equations:
\begin{equation}\label{e426}
\left\{
\aligned
\Pi(x,u,J)=
& \
R_{u^{i_1}u^{i_1}}^{i_1},\\
\Pi(x,u,J)=
& \
R_{u^{i_1}u^{i_1}}^j, \ \ \ \ \ \ \ j\neq i_1,\\
\Pi(x,u,J)=
& \
R_{u^{i_1}u^{i_2}}^{i_1}, \ \ \ \ \ \ \ i_1\neq i_2, \\
\Pi(x,u,J)=
& \
R_{u^{i_1}u^{i_2}}^j, \ \ \ \ \ \ \ i_1\neq i_2, \ \ j\neq i_1, \ \ j\neq i_2.
\endaligned\right.
\end{equation}
This gives the second equation of~(\ref{e415}), completing 
the proof of Lemma~\ref{lem41} and consequently the proof 
of~Theorem~\ref{thm1}. \qed
\endproof

\vfill
\end{document}